\documentclass[12pt]{article}
\usepackage{a4wide,amssymb,amsmath}
\usepackage{amssymb,amsmath}
\usepackage{url}
\usepackage[english]{babel}
\renewcommand{\le}{\leqslant}
\renewcommand{\ge}{\geqslant}
\newcommand{\Supp}{\vv S}
\newcommand{\nz}{\setminus\{0\}}
\newcommand{\ra}{R_{\alpha}}
\newcommand{\vv}[1]{\mathbf{#1}}
\newcommand{\rI}{{\rm I}}
\newcommand{\rH}{{\rm H}}
\newcommand{\cA}{\mathcal{A}}
\newcommand{\cM}{\mathcal{M}}
\newcommand{\cC}{\mathcal{C}}
\newcommand{\cW}{\mathcal{W}}
\newcommand{\RR}{\mathbb{R}}
\newcommand{\ZZ}{\mathbb{Z}}

\newcommand{\NN}{\mathbb{N}}

\newcommand{\R}{\mathbb{R}}
\newcommand{\Z}{\mathbb{Z}}

\newcommand{\N}{\mathbb{N}}

\newcommand{\vx}{\mathbf{x}}
\newcommand{\vy}{\mathbf{y}}
\newcommand{\va}{\mathbf{a}}
\newcommand{\vb}{\mathbf{b}}
\newcommand{\vt}{\mathbf{t}}
\newcommand{\vf}{\mathbf{f}}
\newcommand{\vq}{\mathbf{q}}
\newcommand{\vc}{\mathbf{c}}
\newcommand{\vg}{\mathbf{g}}
\newcommand{\bu}{\mathbf{u}}

\newcommand{\diam}{\mathrm{diam}}
\newcommand{\AAA}{\mathcal{A}}
\newcommand{\HHH}{\mathcal{H}}
\newcommand{\RRR}{\mathcal{R}}
\newcommand{\FFF}{\mathcal{F}}

\newcommand{\proof}{\textit{Proof.}}
\newcommand{\dist}{\mathrm{dist}}
\newcommand{\ve}{\varepsilon}
\newtheorem{theorem}{Theorem}
\newtheorem{lemma}{Lemma}
\newtheorem{definition}{Definition}
\newtheorem{prop}{Proposition}
\newtheorem{corl}{Corollary}
\begin{document}

\title{Inhomogeneous theory of dual Diophantine \\ approximation  on manifolds}

\author{
 Dzmitry Badziahin\footnote{EPSRC PDRA, Grant EP/E061613/1} \and
 Victor Beresnevich\footnote{EPSRC Advanced Research Fellow, grant EP/C54076X/1} \and
 Sanju Velani\footnote{Research partially supported by EPSRC grants EP/E061613/1 and EP/F027028/1 }
}

\date{{\small
{\em Dedicated  to Bob Vaughan on his 65th birthday}}}

\maketitle

\begin{abstract}
The inhomogeneous Groshev type theory for dual Diophantine
approximation on manifolds is developed. In particular,   the notion of   {\em nice}  manifolds is introduced and the  divergence part of the theory is established  for all such manifolds.  Our results naturally incorporate and generalize the homogeneous measure and dimension theorems for non-degenerate manifolds established to date.
The generality of the inhomogeneous aspect considered within  enables us to make a new contribution  even to the classical theory in $\RR^n$.  Furthermore, the  multivariable aspect considered within has natural applications beyond the standard inhomogeneous theory such as to Diophantine problems related to approximation  by algebraic integers.
\end{abstract}

{\small

\noindent\emph{Key words and phrases}: Metric
Diophantine approximation, extremal manifolds,
Groshev type theorem, ubiquitous systems

\medskip

\noindent\emph{AMS Subject classification}: 11J83, 11J13, 11K60

}

\section{Introduction}\label{itro}

\subsection{Motivation and main results \label{mmr}}

Throughout $\RR^+=(0,+\infty)$, $|\cdot|$ denotes the supremum norm,
$\| \cdot\|$ is the distance to the nearest integer and
$\va\cdot\vb:=a_1b_1+\dots+a_nb_n$ is the standard inner product of
vectors $\va=(a_1,\dots,a_n)$ and $\vb=(b_1,\dots,b_n)$ in $\R^n$. Furthermore,
$\Psi:\R^n\to\R^+$ will denote a function such that
\begin{equation}\label{e:001}
    \Psi(a_1,\dots,a_n)\ge\Psi(b_1,\dots,b_n)\quad\text{if
    }|a_i|\le|b_i|\text{ for all }i=1,\dots,n\,,
\end{equation}
and will be referred to as a \emph{multivariable approximating
function}. In the special case when $\Psi(\vv a)=\psi(|\vv a|)$ for
a monotonic function $\psi:\RR^+\to \RR^+$ we simply
refer  to $\psi$  as an \emph{approximating function}.

Given a
multivariable approximating function $\Psi$ and a function $\theta:\R^n\to\R$,  define the set
\begin{equation}\label{e:002}
\cW_n^{\theta}(\Psi)\, :=  \, \left\{\vv y\in\R^n:\begin{array}{l}
\|\va\cdot\vv y+\theta(\vv y)\| < \Psi(\va)\text{ \ \ }\\[0.5ex]
\text{for infinitely many $\vv a\in\Z^n\nz$}
\end{array}\right\}.
\end{equation}
 For obvious reasons points $\vv y$ in $\cW_n^{\theta}(\Psi)$ are referred to as
\emph{$(\Psi,\theta)$-approximable} and when   $\Psi(\vv a)=\psi(|\vv a|)$ we naturally write
$\cW_n^{\theta}(\psi)$ for $\cW_n^{\theta}(\Psi)$.   In the case the function  $\theta$ is constant,
the set $\cW_n^{\theta}(\Psi)$ corresponds to the familiar
inhomogeneous setting within the general theory of dual Diophantine approximation.  In turn,  with $\theta  \equiv 0$  the corresponding set  reduces to the homogeneous  setting and is denoted by
$\cW_n(\Psi)$.    Note that within the homogeneous setting, points in $\R^n$ are approximated by $(n \! - \! 1)$--dimensional rational planes and Groshev's fundamental theorem \cite[\S2.3]{BernikDodson-1999} in the theory of metric Diophantine approximation  provides a beautiful and simple criterion for the `size' of $\cW_n(\psi)$  expressed in terms of $n$-dimensional Lebesgue measure $| \ . \ |_n$. Essentially, for any approximating function $\psi$
$$
\big|\cW_{n}(\psi)\big|_n =
 \left\{ \begin{array}{cl}
 \mbox{\rm
Z{\scriptsize ERO}} &\text{ if}\quad\sum_{t=1}^\infty t^{n-1}\psi(t)
<\infty,  \\[1ex]
\mbox{\rm
F{\scriptsize ULL}} &\text{ if}\quad\sum_{t=1}^\infty t^{n-1}\psi(t)
=\infty  .
\end{array}\right.
$$
Here `\mbox{\rm F{\scriptsize ULL}}' simply means that the complement of the set under consideration is of measure zero. Many years later, and building upon the work of Jarn\'{\i}k,  this criterion was generalized to incorporate  Hausdorff measures  \cite{dettasv}.  For background, precise statements and generalizations to the inhomogeneous and multivariable aspects the reader is refereed to
\cite{Beresnevich-Bernik-Dodson-Velani-Roth, Beresnevich-Dickinson-Velani-06:MR2184760, Beresnevich-Velani-06:MR2264714,Beresnevich-Velani-09-Groshev} and references within.

Let $\cM$ be a manifold in $\R^n$.  In short, our primary goal  is to develop a metric theory for the sets
$\cM   \cap \cW_n^{\theta}(\Psi)$ akin to Groshev's theorem.  The fact that the points $\vv y \in \R^n$ of interest are restricted  to $\cM$  and therefore are of dependent variables, introduces
major difficulties in attempting to describe the measure theoretic
structure of $\cM   \cap \cW_n^{\theta}(\Psi)$ -- even with  $\theta  \equiv 0$ and $\Psi = \psi$.  Trivially, if the dimension  of the manifold is strictly less than $n$ then  $|\cM|_{n}=0$.
Thus, in attempting to develop a Lebesgue theory for $ \cM
\cap \cW_n^{\theta}(\Psi)$  it
is natural to use the induced Lebesgue measure $| \ . \ |_{\cM}$  on
$\cM$.  Trivially, the measure of the complement of $\cM$ with respect to $| \ . \ |_{\cM}$ is zero and so by definition $| \cM  |_{\cM} := \mbox{\rm F{\scriptsize ULL}} $.

In 1998, Kleinbock $\&$ Margulis
\cite{Kleinbock-Margulis-98:MR1652916} established the fundamental
Baker-Sprindzuk conjecture concerning homogeneous Diophantine
approximation on  manifolds.  As a consequence, for non-degenerate manifolds
$ |\cM   \cap \cW_n(\Psi_\epsilon)|_{\cM}= 0  $  where  $\Psi_\epsilon(\vv a) := |\vv a|^{-n - \epsilon}
$ and $\epsilon > 0$. Essentially,
non-degenerate manifolds  are smooth sub-manifolds of $\R^n$ which are sufficiently
curved so as to deviate from any hyperplane. Formally, a  manifold
$\cM$ of dimension $m$ embedded in $\R^n$ is said to be
\emph{non-degenerate}\/ if it arises from a non-degenerate map $\vf:U\to
\R^n$ where $U$ is an open subset of $\R^m$ and $\cM:=\vf(U)$. The
map $\vf:U\to \R^n:\bu\mapsto \vf(\bu)=(f_1(\bu),\dots,f_n(\bu))$ is
said to be \emph{$l$-non-degenerate at $\bu\in U$} if $\vf$  is $l$ times continuously
differentiable on some sufficiently small ball centred at $\bu$
and the partial derivatives of $\vf$ at $\bu$ of orders up to $l$
span $\R^n$. The map $\vf$ is \emph{$l$-non-degenerate} if it is
$l$-non-degenerate at almost every (in terms of $m$-dimensional
Lebesgue measure) point in $U$; in turn the manifold $\cM=\vf(U)$
is also said to be $l$-non-degenerate.  Finally,  we say that  $\vf$ is \emph{non-degenerate} if it is
$l$-non-degenerate for some $l$; in turn the manifold $\cM=\vf(U)$
is also said to be non-degenerate.  It is well known that any real connected analytic
manifold which is not contained in any hyperplane of $\R^n$  is
non-degenerate.

Without a doubt, the proof of  the  Baker-Sprindzuk conjecture has
acted as the catalyst for the subsequent development of the homogeneous
theory of Diophantine approximation on  manifolds. In particular,
the significantly stronger
Groshev type theorem for $\cM   \cap \cW_n(\psi)$ has
been established -- see \cite{Beresnevich-02:MR1905790, Bernik-Kleinbock-Margulis-01:MR1829381} for the zero measure criterion and
\cite{Beresnevich-Bernik-Kleinbock-Margulis-02:MR1944505} for the full measure criterion.  Staying strictly within the homogeneous setting, for recent developments regarding  the `deeper' Hausdorff measure theory and the simultaneous approximation theory  we referee the reader to \cite{Jason, BerS, Beresnevich-Dickinson-Velani-06:MR2184760, Beresnevich-Dickinson-Velani-07:MR2373145, Beresnevich-Velani-NOTE, BerZor, DickinsonDodson-2000a, Vaughan-Velani-2007}  and references within.

Until the recent proof of the inhomogeneous Baker-Sprindzuk
conjecture
\cite{Beresnevich-Velani-Moscow,Beresnevich-Velani-08-Inhom}, the
theory of inhomogeneous Diophantine approximation on manifolds had remained essentially non-existent and ad-hoc -- see \cite{Badziahin-inhom, Bernik-Dickinson-Yuan-1999-inhom, Bernik-Kovalevskaya-06:MR2298911, Ustinov-05:MR2210400,
Ustinov-06:MR2229813}. As a consequence of the measure results in
\cite{Beresnevich-Velani-08-Inhom} we
now know that for any non-degenerate manifold $\cM$  and
$\theta \equiv {\rm constant}$,
\begin{equation}\label{hit}
  |\cM   \cap \cW_n^{\theta}(\Psi_\epsilon)|_{\cM}= 0  \qquad  \forall   \quad  \epsilon > 0  \ .
\end{equation}
Clearly, this statement is far from  the desirable
Groshev type theorem  even for $\cM   \cap \cW_n^{\theta}(\psi) $.  As mentioned above, such a statement exists within
the homogeneous setting. This paper constitutes part of an ongoing  programme
to develop a coherent inhomogeneous theory of Diophantine approximation on manifolds in line with the homogeneous theory.
In the case of  simultaneous approximation on planar curves, the programme has successfully  been carried out in \cite{BobSanVic}.
Here we deal with the dual approximation aspect of the programme.

Our first  result provides a zero Lebesgue measure criterion for $\cM   \cap \cW_n^{\theta}(\Psi) $. It represents the  complete inhomogeneous version of the main result of \cite{Bernik-Kleinbock-Margulis-01:MR1829381}  and it implies  (\ref{hit}) without imposing the condition that  the  `inhomogeneous' function $\theta:\R^n\to\R$ is constant.
Throughout,    $\theta|_\cM$  will denote the
restriction of  the inhomogeneous  function $\theta$ to $\cM$ and as usual, $C^{(n)}$ will denote the set of $n$-times continuously
differentiable functions.

\begin{theorem}\label{t2}
Let $\cM$ be an $l$-non-degenerate manifold in $\R^n$ $(n\ge2)$\/ and
$\theta:\R^n\to\R$ be a function such that $\theta|_\cM  \in C^{(l)}$. Let $\Psi$ be a  multivariable approximating function. Then
\begin{equation*}
\left|\cW_{n}^\theta(\Psi)\cap\cM\right|_\cM=
0\qquad\text{if}\qquad \sum_{\vv a\in\Z^n\nz}\Psi(\vv a)
<\infty.
\end{equation*}
\end{theorem}

For the divergence counterpart, we are able to prove the more general statement in terms of $s$-dimensional Hausdorff measure $\HHH^s$. However, there is a downside in that we impose a `convexity' condition on $\Psi$ which we refer to as property {\bf P}.
For an $n$-tuple $\vv v=(v_1,\dots,v_n)$ of positive numbers satisfying $v_1+\dots+v_n=n$, define the $\vv v$-\emph{quasinorm}\/ $|\,\cdot\,|_{\vv v}$ on $\R^n$ by setting
$$|\vv y|_{\vv v}:=\max_{1\le i\le n}|y_i|^{1/v_i}     \ . $$
{\em A  multivariable approximating function $\Psi$ is said to satisfy property {\bf P} if $\Psi(\vv a)=\psi(|\vv a|_{\vv v})$ for some approximating  function $\psi$ and $\vv v$ as above.}  Trivially, with $\vv v = (1, \ldots, 1)$ we have that $|\vv a|_{\vv v}=|\vv a|$ and  we see that any approximating function $\psi$ satisfies property {\bf P}, where $\psi$ is regarded as the function $\vv a\mapsto\psi(|\vv a|)$.

\begin{theorem}\label{T3}
Let $\cM$ be a non-degenerate manifold in $\R^n$ of dimension
$m$ and let $s > m-1$. Let $\theta:\R^n\to\R$ be a function such that $\theta|_\cM  \in C^{(2)}$ and $\Psi$ be a multivariable approximating function satisfying property {\bf P}. Then
$$
\HHH^s(\cW_n^\theta(\Psi)\cap\cM)=\HHH^s(\cM)\qquad\mbox{if }
\qquad \sum_{\vv a\in\Z^n\nz} |\vv a|  \left(\frac{\Psi(\vv a)}{|\vv a| }\right)^{s+1-m}=\infty.
$$
\end{theorem}

The  above theorem will be derived from a general  statement which significantly broadens the scope of potential applications and is of independent interest. Given a manifold $\cM\subset \R^n$, an $n$-tuple $\vv v=(v_1,\dots,v_n)$ of positive numbers satisfying $v_1+\dots+v_n=n$, $\delta>0$ and $Q>1$, let
$$
\Phi_{\vv v}(Q,\delta):=\Big\{\vy\in \cM: \exists\ \vv a\in\Z^n\setminus\{0\} \text{ such that }\|\va\cdot\vv y\|<\delta Q^{-n}\text{ and } |\vv a|_{\vv v}\le Q\Big\}.
$$
As a consequence of Dirichlet's theorem, $\Phi_{\vv v}(Q,\delta) =  \cM$  if $\delta \geq 1 $. We say that the manifold $\cM$ is \emph{$\vv v$-nice at $\vv y_0\in \cM$}\/ if there is a neighborhood $\Omega\subset \cM$  of $\vy_0$ and constants $0<\delta,\omega<1$ such that for any ball $B\subset \Omega$ we have that
$$
\limsup_{Q\to\infty}\,|\Phi_{\vv v}(Q,\delta)\cap
B|_\cM\le \omega|B|_\cM\,.
$$
The manifold is said to be \emph{$\vv v$-nice}\/ if it is $\vv v$-nice at almost every point in $\cM$. Furthermore, the manifold is said to be \emph{nice}\/ if it is $\vv v$-nice for all choices of $\vv v$.

\begin{theorem}\label{tnice}
Let $\cM$ be a $\vv v$-nice $C^{(2)}$ manifold in $\R^n$\/ of dimension
$m$ and let $s > m-1$. Let $\theta:\R^n\to\R$ be a function such that $\theta|_\cM  \in C^{(2)}$ and $\Psi(\vv a)=\psi(|\vv a|_{\vv v})$ for some approximating function $\psi$. Then
$$
\HHH^s(\cW_n^\theta(\Psi)\cap\cM)=\HHH^s(\cM)\qquad\mbox{if }
\qquad \sum_{\vv a\in\Z^n\nz} |\vv a|  \left(\frac{\Psi(\vv a)}{|\vv a| }\right)^{s+1-m}=\infty.
$$
\end{theorem}

A consequence of Lemma~\ref{lem5} in \S\ref{aux} is that non-degenerate manifolds are nice. Thus
$$
\text{Theorem~\ref{tnice}}\quad\Longrightarrow \quad\text{Theorem~\ref{T3}}.
$$

\subsection{Remarks and Corollaries}

\noindent{\em Remark 1. }
For $s<m$,  the
non-degeneracy of $\cM$ in Theorem \ref{T3}  can be relaxed to the condition that there exists at least one
non-degeneracy point on $\cM$. Also, note that   $\HHH^s(\cM)=\infty$ when $s<m$.

\noindent{\em Remark 2.  } It follows from the definition of Hausdorff measure that
$ \HHH^s(\cW_n^\theta(\psi)\cap\cM)  \leq \HHH^s(\cM) = 0        $
for any $s > m $ irrespective of $\Psi$.  Thus the meat of Theorem~\ref{T3}  is when  $s \leq m$.

\noindent{\em Remark 3. } To the best of our knowledge, Theorem \ref{T3}   with $\cM = \R^n$ and $ \theta \not\equiv {\rm constant}$  is new.  In other words, the theorem  makes a new contribution even to the classical theory of Diophantine approximation of independent variables.

\noindent{\em Remark 4. } We suspect that Property P imposed in the statement of Theorem~\ref{T3} can be safely removed. This constitutes a challenging problem even in the homogeneous case.

\noindent{\em Remark 5. }
Consider the problem of Diophantine approximation on the Veronese curves $\cM:=\{(x,x^2,\dots,x^n):x\in\R\}$, where $n\ge2$. Take
$\theta(x,\dots,x^{n})=x^{n+1}$. Then the inequality in
(\ref{e:002}) becomes
$$
|x^{n+1}+a_{n}x^{n}+\dots+a_1x+a_0| < \Psi(\va)   \,  .
$$
Clearly the function $\theta$ as defined above is $C^{(\infty)}$. In the case when $\Psi(\vv a)=\psi(|\vv a|)$ the corresponding divergence results have been proved by Bugeaud \cite{Bugeaud-?2} and the corresponding convergence results  by Bernik $\&$ Shamukova \cite{Bernik-Shamukova-06:MR2291128, Shamukova-07:MR2397715}. Theorems 1 and 2 naturally extend their results to the   case of multivariable approximating functions $\Psi$.

\bigskip

We now discuss various corollaries of our main theorems which are of independent interest. The following statement is a  direct consequence   of Theorem  \ref{T3} and the fact that any approximating function $\psi$ satisfies property {\bf P}.

\begin{corl}\label{t3}
Let $\cM$ be a non-degenerate manifold in $\R^n$\/ of dimension
$m$ and  $s>m-1$. Let $\theta:\R^n\to\R$ be a function such that $\theta|_\cM  \in C^{(2)}$ and $\psi$ be an  approximating function. Then
$$
\HHH^s(\cW_n^\theta(\psi)\cap\cM)=\HHH^s(\cM)\qquad\mbox{if }
\qquad
\sum_{t=1}^\infty t^n\left(\frac{\psi(t)}{t}\right)^{s+1-m}=\infty.
$$
\end{corl}

\noindent In the case of curves this corollary was first established in \cite{Badziahin-inhom}. In the case $s=m$, the
Hausdorff measure $\HHH^s$ is comparable to the induced  $m$-dimensional Lebesgue measure  $| \ . \ |_{\cM} $ on $\cM$  and Corollary  \ref{t3} represents the  complete inhomogeneous version of the main result of \cite{Beresnevich-Bernik-Kleinbock-Margulis-02:MR1944505}. Furthermore, Theorem \ref{t2} together with Corollary \ref{t3} provides a simple criterion for the `size' of $ \cW_n^\theta(\psi)\cap\cM $  expressed in terms of the induced measure; i.e. the desired inhomogeneous Groshev type theorem for manifolds. More precisely and more generally, under the hypotheses of Theorem \ref{t2} we have that  for any $\Psi$ satisfying property {\bf P}
$$
\big|\cW_{n}^\theta(\Psi)\cap\cM\big|_\cM=
 \left\{ \begin{array}{cl}
0&\text{ if}\quad\sum_{\vv a\in\Z^n\nz}\Psi(\vv a)
<\infty  \\[3ex]
|\cM|_\cM&\text{ if}\quad\sum_{\vv a\in\Z^n\nz}\Psi(\vv a)
=\infty \, .
\end{array}\right.
$$

\noindent In the case  $s < m$, Corollary~\ref{t3}  naturally generalizes the homogeneous result of \cite[Theorem 18]{Beresnevich-Dickinson-Velani-06:MR2184760}.
It is worth stressing that even with $\theta \equiv 0$, Theorem \ref{T3} can be regarded as developing  the  homogenous theory of Diophantine approximation on manifolds by incorporating multivariable approximating functions.

\vspace*{1ex}

Given an approximating function $\psi$, the {\em lower order} $\tau_{\psi}$ of $1/\psi$ is
defined by
$$
\tau_{\psi}:=\liminf_{t\to\infty}\frac{-\log\psi(t)}{\log t}  \;
$$
and indicates the growth of the function $1/\psi$ `near'
infinity. With this definition at hand, it is relatively easy to verify that the divergent sum
condition of Corollary~\ref{t3} is satisfied  whenever $s<m-1+(n+1)/(\tau_{\psi} + 1)$.
It follows from the  definition of Hausdorff measure and dimension that $\dim(\cW_n^\theta(\psi)\cap\cM)\ge s$ if
$\HHH^s(\cW_n^\theta(\psi)\cap\cM)>0$ and $\HHH^s(\cM)> 0 $ if  $s \leq \dim \cM$.
Thus, Corollary~\ref{t3} readily yields the following
 inhomogeneous version of the
dimension result of
\cite{DickinsonDodson-2000a}.

\begin{corl}\label{cor1}
Let $\cM$ be a non-degenerate manifold in $\R^n$\/ of dimension
$m$  and  $\theta:\R^n\to\R$ be a function such that $\theta|_\cM  \in C^{(2)}$. Let $\psi$ be an  approximating function such that
$n  \leq \tau_\psi  < \infty $. Then
\begin{equation}\label{dimdim}
\dim \cW_n^\theta(\psi)\cap\cM \ \ge \ m-1+ \frac{n+1}{\tau_\psi+1} \, .
\end{equation}
\end{corl}

\noindent In the case that  $\theta \equiv {\rm constant} $  and $\psi(t) := t^{-\tau}$ with $\tau>n$,  this dimension statement corresponds to the main result of \cite{Bodyagin-05:MR2212378}. However, Corollary \ref{t3} implies the stronger measure  statement that $\HHH^s(\cW_n^\theta(\psi)\cap\cM)  = \infty $ at $s = m-1+(n+1)/(\tau + 1)$ which in all likelihood is the critical exponent. In a wider context, it would not be unreasonable to expect that the above lower bound for $ \dim \cW_n^\theta(\psi)\cap\cM $ is in fact sharp. Even within the homogenous setting, establishing equality in (\ref{dimdim}) represents a key open problem.  To date the homogeneous problem has been  settled by Bernik \cite{Bernik-1983a} for  Veronese curves  and by R.C.\,Baker \cite{Baker-1978} for non-degenerate planar curves.  For non-degenerate curves in $\R^n$ the current results are limited to situation that $\tau_\psi \leq n+\frac1{4n}$ -- see \cite{Beresnevich-Bernik-Dodson-02:MR2069553}. Most recently, the  inhomogeneous version of Baker's result  has  been established  in \cite{Badziahin-inhom}.  In other words, if $\cM$ is a non-degenerate planar curve then in (\ref{dimdim})  we have equality.

\subsection{Possible developments}

\noindent \emph{Affine subspaces.} The homogeneous Groshev type theorems are established for lines passing through the origin in \cite{Beresnevich-Bernik-Dickinson-Dodson-00:MR1820985}
and for hyperplanes in \cite{Ghosh-05:MR2156655}. It is likely that the techniques developed in this paper can be used to extend these results to the inhomogeneous setting. Note that affine subspaces are degenerate everywhere and so Theorems~\ref{t2} $\&$ \ref{T3} are not applicable. It is worth mentioning that the general case, in which the dimension of planes is arbitrary, is unresolved even in the homogeneous setting -- see \cite{Kleinbock-03:MR1982150} for further details.

\bigskip

\noindent \emph{Manifolds in affine subspaces.} By definition, any manifold contained in a proper affine subspace of $\R^n$ is degenerate. Nevertheless the `extremal' theory of homogeneous Diophantine approximation for such manifolds has been developed in \cite{Kleinbock-03:MR1982150}. A natural problem is to develop the analogous inhomogeneous theory and more generally obtain a Groshev type theory. For the same reason as above,  Theorems~\ref{t2} $\&$ \ref{T3} are not applicable in this setting.

\bigskip

\noindent \emph{The $p$-adic setting.}
The homogeneous Groshev type theorems have recently been established in \cite{Mohammadi-Golsefidy-09, Mohammadi-Golsefidy-preprint} for the `$S$-arithmetic' setting. This builds upon the `extremality' results of Kleinbock and Tomanov \cite{Kleinbock-Tomanov-07:MR2314053} and includes the more familiar $p$-adic case. In all likelihood the techniques developed in this paper can be used to extend the homogeneous $S$-arithmetic results to the inhomogeneous setting.    For  inhomogeneous  $p$-adic results restricted to  Veronese curves see \cite{Bernik-Dickinson-Yuan-1999-inhom, Budarina-Zorin-09:MR2530194, Ustinov-05:MR2210400}.

\bigskip

\noindent\emph{The non-monotonic setting.}  By definition, any approximating function $\psi$ is  monotonic.  Thus, monotonicity is implicitly assumed within the context of the classical Groshev theorem as stated in \S\ref{mmr}.  Recently in \cite{Beresnevich-Velani-09-Groshev},   this classical result has been freed from all unnecessary monotonicity constraints.  Naturally, it would be highly desirable to obtain analogous statements for Diophantine approximation on manifolds. This in full generality is a difficult problem. Even in the case $\Psi(\vv a)=\psi(|\vv a|)$,  to remove the implicit monotonicity assumption from Theorems~\ref{t2} $\&$ \ref{T3}  is believed to be currently out of reach.  For homogeneous convergent Groshev type results without monotonicity but restricted to non-degenerate curves in $\R^n$ see \cite{Beresnevich-05:MR2110504, Budarina-Dickinson-09:MR2481999}. In the first instance it would be interesting to extend these homogeneous results for curves  to the inhomogeneous setting.

\subsection{Global assumptions and useful conventions}\label{useful}

In the course of proving our results we will conveniently and without loss of generality  assume that the
manifold $\cM$ under  consideration is immersed in $\R^n$ via a smooth map $\vv f=(f_1,\dots,f_n):U\to\R^n$ defined on a ball  $U\subset\R^m$.
Thus, $\cM=\big\{\vv f(\vx):\vx\in U\big\}$.
Furthermore, in view of the Implicit Function Theorem we can assume that
$$f_i(\vv x)=x_i \qquad {\rm for}  \quad  i=1,\dots,m  \, . $$
In other words,  $\vv f$ is a Monge parameterisation of $\cM$.   Note that this implies that  $\vv f$ is locally bi-Lipschitz.

Let $\cA_{\vv f}(\Psi,\theta)$ denote the projection of $\cW_n^\theta(\Psi)\cap\cM$ onto $U$; that is
$$
\cA_{\vv f}(\Psi,\theta):=\big\{\vx\in U\,:\,\vf(\vx)\in \cW_n^\theta(\Psi)\big\}\,.
$$
Thus,  a point  $\vv x\in \cA_{\vv f}(\Psi,\theta)$ if and only if the point  $\vv f(\vv x)\in\cM$ is $(\Psi,\theta_{\vv f}(\vv x))$-approximable with $\theta_{\vv f}(\vv x):=\theta(\vv f(\vv x))$. For convenience and clarity we will drop the subscript from $\theta_{\vv f}$.
In the case when $\Psi(\vv a)=\psi(|\vv a|)$ for some
approximating function $\psi$ we write $\cA_{\vv f}(\psi,\theta)$
for $\cA_{\vv f}(\Psi,\theta)$. A consequence of the fact that $\vv f$ is locally bi-Lipschitz is that Theorems~\ref{t2}--\ref{tnice} can be equally stated in terms of $\cA_{\vv
f}(\Psi,\theta)$.  Indeed the proof of the theorems will make use of this alternative formulation.

In the case of Theorem~\ref{t2} the functions $\vv f$ and $\theta$ are $C^{(l)}$. Thus we  can assume without loss of generality that there is a constant $C_0>0$ depending only on $U$, $\vv f$ and $\theta$  such that
\begin{equation}\label{bounds}
    \max_{0\le i\le l}\,\,\sup_{\vv x\in U}|\vv f^{(i)}(\vv x)|\le C_0\qquad \text{and}\qquad\max_{0\le i\le l}\,\,\sup_{\vv x\in U}|\theta^{(i)}(\vv x)|\le C_0.
\end{equation}
In the case of Theorems~\ref{T3} $\&$  \ref{tnice} the functions $\vv f$ and $\theta$ are $C^{(2)}$ and therefore  without loss of generality we can assume (\ref{bounds}) with $l=2$.

\bigskip

\noindent{\em Notation.} The Vinogradov symbols $\ll$
and $\gg$ will be used to indicate an inequality with an
unspecified positive multiplicative constant. If $a \ll b$ and $a
\gg b$ we write $ a \asymp b $, and  say that the quantities $a$
and $b$ are comparable.  We denote by  $B=B(\vx,r)$ the ball centred at  $\vx\in\R^m$ with radius $r$.  For any real number $\lambda>0$, we let  $\lambda
 B$  denote the ball $B$ scaled by a factor $\lambda$; i.e.  $\lambda B(\vx,
 r):= B(\vx, \lambda r)$.

\section{The convergence theory}\label{convergence}

The goal is to prove  Theorem \ref{t2}.  Thus,
throughout  $\Psi$ is a multivariable approximating function satisfying the convergent sum condition
\begin{equation}\label{e:017}
\sum_{\vv a\in\Z^n\nz}\Psi(\vv a)
<\infty \, .
\end{equation}
In view of the discussion of \S\ref{useful} the goal is equivalent to establishing $|\cA_\vf(\Psi,\theta)|_m=0$. Note that the set $\cA_\vf(\Psi,\theta)$  can be written as
\begin{equation*}
\cA_\vf(\Psi,\theta)=\limsup_{|\vv
a|\to\infty}A_\vf(\va,\Psi,\theta):=\bigcap_{h=1}^\infty \
\bigcup_{|\va|\ge h} A_\vf(\va,\Psi,\theta),
\end{equation*}
where $$A_\vf(\va,\Psi,\theta):=\big\{\vv x\in U :
\|\va\cdot\vv f(\vv x)+\theta(\vv x)\| < \Psi(\va)
\big\}   \, . $$

\noindent For each $\vv a\in\Z^n\nz$ it is convenient to decompose  the set $A_\vf(\va,\Psi,\theta)$  into the following two subsets
\begin{equation}\label{e:020}
A^1_\vf(\va,\Psi,\theta):=\big\{\vv x\in A(\va,\Psi,\theta)\;:\; |\nabla
(\vf\cdot \va+\theta)(\vx)|\ge C_1\,|\va|^{1/2}\,\big\}
\end{equation}
and
\begin{equation*}
A^2_\vf(\va,\Psi,\theta):=\big\{\vv x\in A(\va,\Psi,\theta)\;:\; |\nabla
(\vf\cdot\va+\theta)(\vx)|< C_1\,|\va|^{1/2}\,\big\}\, .
\end{equation*}
Here $\nabla$  as usual  denotes the gradient  operator and
\begin{equation}\label{C_1}
C_1:=\sqrt{(n+1)mC_0}
\end{equation}
where  $C_0$ is as  in (\ref{bounds}). Obviously
$$
\cA_\vf(\Psi,\theta)=\cA_\vf^{1}(\Psi,\theta)\cup
\cA_\vf^{2}(\Psi,\theta),
$$
where
$$
\cA_\vf^{i}(\Psi,\theta)=\limsup_{|\vv
a|\to\infty}A^i_\vf(\va,\Psi,\theta):=\bigcap_{h=1}^\infty \
\bigcup_{|\va|\ge h} A^i_\vf(\va,\Psi,\theta) \qquad   (i=1,2)\,.
$$

\noindent The desired statement that $|\cA_\vf(\Psi,\theta)|_m=0$ will follow by establishing  the  separate cases:\\[3ex]
\hspace*{10.5ex}\textbf{Case A}  \ \ \ \ $|\cA^1_\vf(\Psi,\theta)|_m=0$ \ \ \ \\[3ex]
\hspace*{10.5ex}\textbf{Case B} \ \ \ \ $|\cA^2_\vf(\Psi,\theta)|_m=0$.

\bigskip

\subsection{Establishing Case A}\label{case_i}

The aim is to show that $|\cA_\vf^1(\Psi,\theta)|_m=0$.
This will follow as a consequence  of Theorem~1.3 from
\cite{Bernik-Kleinbock-Margulis-01:MR1829381} which is now explicitly stated using slightly different notation.

\begin{theorem}[Bernik, Kleinbock \& Margulis]\label{bkm_th}
Let $B\subset\RR^m$ be a ball of radius $r>0$ and let
$\vg= (g_1,g_2,\ldots, g_{n+1})\in C^{(2)}(2B)$. Fix $\delta>0$ and suppose that
\begin{equation}\label{e:022}
L:=\max_{1\le i,j\le m}\sup_{\vv x\in 2B} \left|\frac{\partial^2 \vg(\vv x)}{\partial x_i\partial x_j}\right|\ < \ \infty\,.
\end{equation}
Then for every $\vq\in \ZZ^{n+1}$ such that
\begin{equation}\label{e:023}
|\vq|\ge \frac{1}{4(n+1)Lr^2}
\end{equation}
the set of $\vv x\in B$ satisfying the system of inequalities
\begin{equation}\label{e:024}
\left\{\begin{array}{ccl} \|\vg(\vv x)\cdot \vq\|&<&\delta \\[2ex]
|\nabla \vg(\vv x)\cdot\vq|&\ge&
\big((n+1)mL\,|\vq|\big)^{1/2}
\end{array}\right.
\end{equation}
has measure at most $K  \, \delta |B|_m$, where $K$ is a constant
depending only on $m$.
\end{theorem}

\bigskip

With the above theorem at our disposal, consider any non-empty open ball $B$ such that $2B\subset U$. Let $\vg= (f_1,f_2,\ldots,f_n,\theta)$ and
$\vq= (a_1,\ldots a_n, 1)$ where $\vv a=(a_1,\dots,a_n)\in\Z^n\nz$. Then,  in view of (\ref{bounds}), we have that (\ref{e:022}) is automatically satisfied. Furthermore, (\ref{e:023}) holds
for all except finitely many $\va\in\Z^n\nz$. In view of (\ref{bounds}) and (\ref{C_1}), the lower bound inequality of (\ref{e:024}) is implied by the inequality associated with (\ref{e:020}). Therefore, $A_\vf^1(\va,\Psi,\theta)\cap B$ is contained in the set defined by (\ref{e:024}) with $\delta := \Psi(\va)$. It now follows via  Theorem~\ref{bkm_th}, that
$$
|A_\vf^1(\va,\Psi,\theta)\cap B|_m\ll \Psi(\va)  \
$$
where the implied constant is independent of $\vv a$. This together with  (\ref{e:017}) and  the Borel-Cantelli lemma readily implies that
$|\cA_\vf^1(\Psi,\theta)\cap B|_m=0$. Now simply observe that the open balls $B$ such that $2B\subset U$ cover the whole of $U$. The upshot is that  $|\cA_\vf^1(\Psi,\theta)|_m=0$ as required.

\subsection{Preliminaries for establishing Case B}\label{prelimB}

Establishing Case B relies upon the recent transference technique introduced in \cite{Beresnevich-Velani-08-Inhom}  and the properties of $(C,\alpha)$-good functions introduced by  Kleinbock $\&$ Margulis in \cite{Kleinbock-Margulis-98:MR1652916}.

\subsubsection{Good functions}\label{Calphagood}
 The following formal definition can be found in \cite{Kleinbock-Margulis-98:MR1652916}.

\begin{definition}\label{def2}\rm
Let $C$ and $\alpha$ be positive numbers and $f:V\to\R$ be a
function defined on an open subset $V$ of $\;\RR^m$. Then $f$ is
called {\em $(C,\alpha)$-good on $V$} if for any open ball $B\subset V$
and any $\varepsilon>0$ one has that
\begin{equation}\label{eq3}
\left|\Big\{x\in B\;:\; |f(x)| \, <\, \varepsilon\;  \sup_{x\in
B}|f(x)|\Big\}\right|_m \ \le  \  C\varepsilon^\alpha |B|_m.
\end{equation}
\end{definition}

\medskip

We now recall various useful properties of $(C,\alpha)$-good functions.

\medskip

\begin{lemma}{\bf  (\cite[Lemma~3.1]{Bernik-Kleinbock-Margulis-01:MR1829381})}
\label{BKM_lem3_1} \
\begin{itemize}
\item[\rm{(a)}] If $f$ is $(C,\alpha)$-good  on $V$ then so is
$\gamma f$ for any $\gamma\in \RR$.
\item[\rm{(b)}] If $f$ and $g$  are $(C,\alpha)$-good on $V$ then so is
$\max \{|f|, |g|\}$.
\item[\rm{(c)}] If $f$ is $(C,\alpha)$-good on $V$ then $f$ is
$(C',\alpha')$-good on $V'$   for every $C'\ge C, \alpha'\le \alpha$
and $V'\subseteq V$.
\item[{\rm(d)}]If $f$ is $(C, \alpha)$-good on $V$ and $c_1 \leq \frac{|f(x)|}{|g(x)|} \leq c_2$ for all $x \in V$,  then $g$ is $\left(C(c_2/c_1)^{\alpha}, \alpha\right)$-good on $V$.
\end{itemize}
\end{lemma}

\noindent The next lemma is the key tool for establishing that a given function is $(C,\alpha)$-good. The following notation is needed to state the lemma. An $m$-tuple $\beta=(\beta_1,\dots,\beta_m)$ of non-negative integers will be referred to as a multiindex and we let $|\beta|_*:=\beta_1+\dots+\beta_m$.  Given a multiindex $\beta$, let
 $$\partial_\beta := \dfrac{\partial^{|\beta|_*}}{\partial x_m^{\beta_1}\cdots\partial x_m^{\beta_m}}   \quad { \rm and } \quad \partial_i^k:=\dfrac{\partial^k}{\partial x_i^k}  .  $$

\begin{lemma}{\bf  (\cite[Lemma~3.3]{Bernik-Kleinbock-Margulis-01:MR1829381})}\label{BKM_lem3_3}
Let $U$ be an open subset of\/ $\RR^m$ and let $g\in C^{(k)}(U)$ be such
that for some constants $A_1,A_2>0$
\begin{equation}\label{e:026}
|\partial_\beta g(\vx)|\le A_1\quad \forall\;\beta\mbox{ with
}|\beta|_*\le k,
\end{equation}
and
\begin{equation}\label{e:027}
|\partial_i^k g(\vx)|\ge A_2\quad \forall\; i=1,\ldots,m
\end{equation}
for all $\vx\in U$. Also let $V$ be a subset of\/ $U$ such that
whenever a ball $B$ lies in $V$ any cube circumscribed around $B$
is contained in $U$. Then $g$ is $(C,\frac{1}{mk})$-good on $V$ for
some explicit positive constant $C$ depending on $A_1$, $A_2$, $m$ and $k$ only.
\end{lemma}

The following proposition\footnote{In Proposition~\ref{prop2} we assume that $\FFF$ is compact. This assumption is not made in Proposition~3.4 of \cite{Bernik-Kleinbock-Margulis-01:MR1829381} although it is used in its proof. Note that the compactness of $\FFF$ does not follow from the assumption that $\{\nabla f:f\in\FFF\}$ is compact. In fact, the family $\FFF$ defined in Corollary~3.5 of \cite{Bernik-Kleinbock-Margulis-01:MR1829381}, which is the main application of \cite[Proposition~3.4]{Bernik-Kleinbock-Margulis-01:MR1829381}, is not compact. The proof of the corollary as given in \cite{Bernik-Kleinbock-Margulis-01:MR1829381} is therefore incomplete. Nevertheless, the corollary as stated is correct. These issues are carried over unaddressed into Theorem~4.5 of \cite{Mohammadi-Golsefidy-09}. In this paper the issues are addressed by our Proposition~\ref{prop2} and Corollary~\ref{BKM_corl3_5}.} is a generalization of Proposition~3.4 from \cite{Bernik-Kleinbock-Margulis-01:MR1829381}.

\begin{prop}\label{prop2}
Let $U$ be an open subset of $\RR^m$, $\vx_0\in U$ and  let
$\FFF\subset C^{(l)}(U)$ be a compact family of functions $f:U\to \RR$ for some  $l \ge 2$ . Assume also that
\begin{equation}\label{e:029}
\inf_{f\in\FFF}\, \max_{0<|\beta|_*\le l}|\,\partial_\beta f(\vx_0)|>0\,.
\end{equation}
Then there exists a neighborhood $V\subset U$ of $\vx_0$ and
positive constants $C$ and $\delta$ satisfying the following property. For any $\Theta\in C^{(l)}(U)$ such that
\begin{equation}\label{e:030}
\sup_{\vv x\in U}\, \max_{|\beta|_*\le l}\,|\partial_\beta
\Theta(\vx)| \, \le \,  \delta
\end{equation}
and any $f\in\FFF$ we have that
\begin{itemize}
\item[\rm{(a)}] $f+\Theta$ is $\big(C,\frac{1}{ml}\big)$-good on $V$, \\
\item[\rm{(b)}] $|\nabla (f+\Theta)|$ is $\big(C,\frac{1}{m(l-1)}\big)$-good
on $V$.
\end{itemize}
\end{prop}

\proof{} The proof is a modification of the ideas used to establish Proposition~3.4 in \cite{Bernik-Kleinbock-Margulis-01:MR1829381}.  First of all note that in view of (\ref{e:029}), there exists a constant $C_1>0$ such that for any
$f\in\FFF$ one can find a multiindex $\beta$ with $0<|\beta|_*=k\le l$, where $k=k(f)$, such that
\begin{equation}\label{beta}
|\partial_\beta f(\vx_0)|\ge C_1.
\end{equation}
Since the number of different $\beta$'s is finite, without loss of generality we can assume that $\beta$ appearing in (\ref{beta}) is the same for all $f\in\FFF$. By an appropriate rotation of the coordinate system one can ensure that
\begin{equation}\label{ineq_ikf}
|\tilde\partial_i^k f(\vx_0)|\ge C_2
\end{equation}
for all $i=1,\ldots,m$ and some positive $C_2$ independent of $f$. Here $\tilde\partial$ denotes differentiation with respect to the rotated coordinate system. Also, by (\ref{e:030}) there exists a constant $c=c(l)>1$ such that
\begin{equation}\label{e:030+}
\sup_{\vv x\in U}\, \max_{|\beta|_*\le l}\,|\tilde\partial_\beta
\Theta(\vx)|\le c\delta.
\end{equation}
Now take $\delta:= C_2/(2c)$. Then, by (\ref{ineq_ikf}) and (\ref{e:030+}), for any $f\in\FFF$ we have that
$$
|\tilde\partial_i^k (f+\Theta)(\vx_0)|\ge \delta \quad\mbox{ for all
}i=1,\ldots,m.
$$
Then, by the continuity of derivatives of $f+\Theta$ and the compactness of $\FFF$, we can choose a neighborhood $V'\subset U$ of $\vv x_0$ and positive constants $A_1,A_2$ independent of $f$ such that (\ref{e:026}) and (\ref{e:027}) with $\partial$ replaced by ${\tilde\partial}$ hold for all $\vx\in V'$ and all $g=f+\Theta$. Finally, let $V$ be a smaller neighborhood of $\vx_0$ such that whenever a ball $B$ lies in $V$, the cube $\widetilde{B}$ circumscribed around $B$ is contained in $V'$. Then,
 on applying Lemma~\ref{BKM_lem3_3} establishes  part~(a) of Proposition~\ref{prop2}.

 Regarding part (b), first assume that $k$ appearing in (\ref{ineq_ikf}) is at least $2$. Since $\FFF$ is compact and differentiation is a continuous map from $C^{(l)}(U)$ to $C^{(l-1)}(U)$, we have that for every $i=1,\dots,m$
\begin{equation*}
\FFF_i:=\big\{\tilde\partial_if\;:\; f\in\FFF\big\}\quad\mbox{ is compact in }C^{(l-1)}(U).
\end{equation*}
In view of the definition of $\FFF$ condition (\ref{e:029}) holds when $l$ is replaced by $l-1$ and $\FFF$ is replaced by $\FFF_i$.
Therefore, the arguments used to prove  part (a) apply to $\FFF_i$ and we conclude  that for every $f\in\FFF_i$ the function $\tilde\partial_i (f+\Theta)$ is $\big(C_i, \frac{1}{m(l-1)}\big)$-good
on some neighborhood $V_i$ of $\vv x_0$. It follows via Lemma~\ref{BKM_lem3_1}, that $|\tilde\nabla(f+\Theta)|$ is $\big(\tilde C,\frac{1}{m(l-1)}\big)$-good with $\tilde C=\max_i C_i$, $V=\cap_i V_i$ and $f\in\FFF$. Naturally,  $\tilde\nabla$ denotes the  gradient operator with respect to the rotated coordinate system. Now simply notice  that the quantity  $$\frac{|\nabla (f+\Theta)(\vx)|}{ |\tilde\nabla(f+\Theta)(\vx)|}$$ for all $\vx \in V$ is bounded between two positive constants.
Hence, by making use of part (d) of  Lemma~\ref{BKM_lem3_1} we obtain the statement of part (b) of  Proposition~\ref{prop2}.

\noindent It remains to consider the case when $k$ appearing in \eqref{ineq_ikf} is
equal to  $1$. Let $A_1,A_2$ and $V$ be defined as in the proof of part (a) above. Then,
\begin{equation}\label{eq1}
A_2\le |\tilde\nabla (f+\Theta)(\vx)|\le A_1\ \ \ \text{ for all $\vx\in V$}.
\end{equation}
In view of part (d) of Lemma~\ref{BKM_lem3_1} and the definition of $(C,\alpha)$-good functions, to complete the proof it suffices to verify that
\begin{equation}\label{vb10}
\left| \{\vx\in B\;:\; |\tilde\nabla
(f(\vx)+\Theta(\vx))|<\varepsilon  \; \sup_{\vy\in B}|\tilde\nabla
(f(\vy)+\Theta(\vy))|\}\right|_m \le
\left(\frac{A_1}{A_2}\right)^{\frac{1}{l-1}}\varepsilon^{\frac{1}{l-1}} \, |B|_m
\end{equation}
for any positive $\varepsilon$ and any $B\subset V$. Firstly, note that if $\varepsilon\ge A_2 / A_1$ then
the r.h.s. of (\ref{vb10}) is at least $|B|_m$ and so (\ref{vb10}) is obviously true. Thus, suppose that  $\varepsilon< A_2 / A_1 $,  Then in view of  (\ref{eq1}), the set on the l.h.s. of (\ref{vb10}) is empty and again (\ref{vb10}) is trivially satisfied.   This thereby completes the proof of the proposition.
\hfill  $\Box$

\bigskip

\begin{corl}\label{BKM_corl3_5}
Let\, $U$ be an open subset of\/ $\RR^m$, $\vx_0\in U$ and let\/
$\vf=(f_1,\ldots,f_n):U\to \RR^n$ be $l$-nondegenerate at
$\vx_0$ for some  $l\ge2$. Let $\theta\in C^{(l)}(U)$. Then there exists a neighborhood $V\subset U$ of\/ $\vx_0$ and positive constants $C$ and $H_0$ such that for any\/ $\va\in\RR^n$ satisfying $|\va|\ge H_0$
\begin{itemize}
\item[$(a)$] $a_0+\va\cdot \vf+\theta$ is
$(C,\frac{1}{ml})$-good on $V$ for every $a_0 \in \RR$, and
\vspace*{1ex}

\item[$(b)$] $|\nabla (\va\cdot \vf+\theta)|$ is
$(C,\frac{1}{m(l-1)})$-good on $V$.
\end{itemize}
\end{corl}

\proof{} To start with choose the neighborhood $V\subset U$ of $\vv x_0$ so that $\vv f$ and $\theta$ are bounded on $V$. Then there exists  a positive constant $K$ such that
\begin{equation}\label{eq2}
\sup_{\vv x\in V}|\vv f(\vv x)|\le K/(n+1)\qquad\text{and}\qquad
\sup_{\vv x\in V}|\theta(\vv x)|\le K/(n+1).
\end{equation}
Let $f$ be the function given by $f(\vv x):=a_0+\va\cdot \vf(\vx)+\theta(\vx)$. Assume for the moment that  $|a_0|\ge2K|\vv a|$.  Then, on using (\ref{eq2}) we find  that
\begin{equation*}
    \sup_{\vv x\in B}|f(\vv x)| \ \le \ 3\inf_{\vv x\in B}|f(\vv x)|
\end{equation*}
for any ball $B\subset V$.
Therefore, if $\varepsilon< 1/3$ then the set on the l.h.s. of (\ref{eq3}) is empty and (\ref{eq3}) is trivially satisfied with any positive $C$ and $\alpha$. On the other hand, if $\varepsilon\ge 1/3$, then (\ref{eq3})
is obviously true for any $C\ge 3$ and any positive $\alpha\le 1$. The upshot is that  part (a) of the corollary holds for any $C\ge 3$ and $0<\alpha\le 1$ whenever $|a_0|\ge 2K|\vv a|$. Thus, without loss of generality we will assume that $|a_0|\le 2K|\vv a|$.

Let $\FFF$ be the collection of functions of the form
$\vc \cdot\vf(\vx)+c_0$, where $\vv c\in\R^n$ such that  $|\vc|=1$ and $|c_0|\le 2K$. Using the compactness of the set
$$
\{\vv c\in\R^n:|\vv c|=1\}\times\{c_0\in\R:|c_0|\le 2K\}  \, ,
$$
one readily verifies that $\FFF$ is compact in $C^{(l)}(U)$.  This together with the fact that $\vf$ is non-degenerate at $\vx_0$   ensures that $\FFF$ satisfies (\ref{e:029}).   Next note that by shrinking the neighborhood $V$ of $\vv x_0$ if necessary, we have  that
$$
\sup_{\vx\in V}\max_{|\beta|_*\le l}|\partial_\beta \theta(\vx)|\le M
$$
for some positive constant $M$.
Now let $C$ and $\delta$ be the contants associated with  Proposition~\ref{prop2}  and let
 $$H_0 \, :=  \, M/\delta  \, .$$
Consider an arbitrary vector $\va\in\RR^n$ with $|\va|\ge H_0$ and any real number  $a_0$ such that $|a_0|\le 2K|\vv a|$. Then, $\Theta$ given by $\Theta(\vv x):= \theta(\vx)/|\vv a|$ satisfies
(\ref{e:030}) and $$f:  \vv x \to f(\vv x):=|\vv a|^{-1}(a_0+\vv f(\vv x)\cdot\vv a)$$ belongs to the compact family $\FFF$. In view of  Proposition~\ref{prop2}, the function $f+\Theta$ given by $f(\vv x)  +\Theta(\vv x)  =|\vv a|^{-1}(a_0+\vv f(\vv x)\cdot\vv a+\theta(\vv x))$ satisfies the desired conclusions of the corollary.  The assertions for the function without the $|\vv a|^{-1}$ multiplier are a  simple consequence of part (a) of Lemma~\ref{BKM_lem3_1}.
\hfill  $\Box$

\begin{prop}\label{prop2+}
Let $U$, $\vx_0$ and $\FFF$ be as in Proposition~\ref{prop2} and suppose that \eqref{e:029} is valid. Then for any neighborhood $V\subset U$ of $\vx_0$,  we have   that  $$\inf_{f\in\FFF}\sup_{\vv x\in V}|f(\vv x)|>0 \, . $$
\end{prop}

\proof{} In view of  (\ref{e:029}) it follows that $\|f\|_V:=\sup_{\vv x\in V}|f(\vv x)|>0$ for every $f\in\FFF$ and any neighborhood $V\subset U$ of $\vx_0$. The map
$f\mapsto\|f\|_V$  is continuous with respect to the $C^{(0)}$ norm. By the compactness of $\FFF$, we have that $\inf_{f\in\FFF}\|f\|_V=\|f_0\|_V$ for some $f_0\in\FFF$. The  claim of the proposition now  follows on combining these facts.
\hfill  $\Box$

\begin{corl}\label{BKM_corl3_5+}
Let\, $U$, $\vx_0$, $\vf$ and $\theta$ be as in Corollary~\ref{BKM_corl3_5}. Then for every sufficiently small neighborhood $V\subset U$ of\/ $\vx_0$,  there exists  $H_0>1$ such that
$$
\inf_{\stackrel{\scriptstyle (\vv a,a_0)\in\R^{n+1}}{|\vv a|\ge H_0}} \ \sup_{\vv x\in V}\ |a_0+\va\cdot \vf(\vx)+\theta(\vx)|>0.
$$
\end{corl}

\proof{} Consider any neighborhood  $V\subset U$ of\/ $\vx_0$ for which the  inequalities given by \eqref{eq2} are satisfied for some $K>0$. Let $f$ denote the function given by $f(\vv x):=a_0+\va\cdot \vf(\vx)+\theta(\vx)$. Notice that if $|a_0|\ge2K|\vv a|$, then in view of   (\ref{eq2}) it follows  that
$$\sup_{\vv x\in V}|f(\vv x)| \ \ge \ K\, H_0>K>0$$
for any $(\vv a,a_0)\in\R^{n+1}$ with $|\vv a|\ge H_0 > 1$ and $|a_0|\ge2K|\vv a|$.  Thus for the rest of the proof we may assume that $|a_0|\le 2K|\vv a|$.

\noindent As in the proof of Corollary \ref{BKM_corl3_5}, let $\FFF$ be the collection of functions of the form
$\vc \cdot\vf(\vx)+c_0$, where $\vv c\in\R^n$ such that  $|\vc|=1$ and $|c_0|\le 2K$.  Then $\FFF$ is a compact subset of $C^{(l)}(U)$  and  since $\vf$ is non-degenerate at $\vx_0$, we have that  $\FFF$ satisfies (\ref{e:029}).
Thus,  Proposition~\ref{prop2+} implies that
$M:=\inf_{f\in\FFF}\ \sup_{\vv x\in V}\ |f(\vx)|>0$.
Therefore, for any $(\vv a,a_0)\in\R^{n+1}$ with $|\vv a|\ge H_0 > 1$ and $|a_0|\le 2K|\vv a|$ we have that
\begin{equation}\label{vb100}
\sup_{\vv x\in V}\ |a_0+\va\cdot \vf(\vx)|\ge M H_0.
\end{equation}
 Now take $H_0  >  \max \{ 1,  K/M \}$. Then, by (\ref{eq2}) and (\ref{vb100}) it follows that
$$\sup_{\vv x\in V}\ |a_0+\va\cdot \vf(\vx)+\theta(\vv x)|\ge M H_0/2$$ and this completes the  proof of the corollary.
\hfill  $\Box$

\subsubsection{Inhomogeneous Transference Principle}\label{ITP}

In this section we describe a simplified version of the  Inhomogeneous Transference Principle introduced in  \cite[Section 5]{Beresnevich-Velani-08-Inhom}. The simplified version  takes into consideration the specific applications that we have in mind.  Throughout,  $V$ denotes  a finite open ball in $\R^m$ and  $\mu$ is $m$-dimensional Lebesgue measure restricted to $V$. Clearly the support of $\mu$ is the closure $\overline V$ of $V$.  For consistency with the notation used in  \cite{Beresnevich-Velani-08-Inhom},   will be write  $\Supp$ for $\overline V$.

Let $\vv T$ and $\AAA$ be two countable `indexing' sets and let $\rH$ and $\rI$ be two maps from $\vv T\times \AAA\times \RR^+$ into the set of open subsets of $\R^m$ such that
\begin{equation}\label{HI}
\rH\;:\; (\vt,\alpha, \varepsilon)\mapsto \rH_\vt(\alpha,\varepsilon)
\qquad\text{and}\qquad
\rI\;:\; (\vt,\alpha, \varepsilon)\mapsto \rI_\vt(\alpha,\varepsilon).
\end{equation}
Let $\Phi$ denote a set of functions $\phi\,:\,\vv T\to
\RR^+$. For $\phi\in \Phi$, consider the $\limsup$ sets

\begin{equation}\label{vb025}
\Lambda_\rI(\phi):=\limsup\limits_{\vt\in\vv
T}\bigcup\limits_{\alpha\in\AAA}\rI_\vt(\alpha,\phi(\vt))\qquad\text{and}\qquad
\Lambda_\rH(\phi):=\displaystyle\limsup\limits_{\vt\in\vv T}
\bigcup\limits_{\alpha\in\AAA}\rH_\vt(\alpha,\phi(\vt)).
\end{equation}

\medskip

\noindent The following  two key properties enables us to transfer zero $\mu$-measure statements for the `homogenous' $\limsup$ sets $\Lambda_\rH(\phi)$ to the  `inhomogenous' $\limsup$ sets $\Lambda_\rI(\phi)$.

\medskip

\noindent{\bf Intersection Property:} The triple $(\rH,\rI,\Phi)$ is said to satisfy the intersection property if for any $\phi\in \Phi$ there exists $\phi^*\in\Phi$ such that for all but finitely many $\vt\in\vv T$ and all distinct $\alpha, \alpha'\in \AAA$
\begin{equation}\label{e:033}
\rI_\vt(\alpha,\phi(\vt))\cap \rI_\vt(\alpha',\phi(\vt))\subset
\bigcup_{\alpha''\in \AAA} \rH_\vt(\alpha'',\phi^*(\vt)).
\end{equation}

\medskip

\noindent{\bf Contracting Property:}
We say that $\mu$ is contracting with respect to $(\rI,\Phi)$ if for any $\phi\in
\Phi$ there exists $\phi^+\in \Phi$ and a sequence of positive
numbers $\{k_\vt\}_{\vt\in\vv T}$ such that
\begin{equation}\label{k_t}
\sum_{\vt\in\vv T} k_\vt<\infty
\end{equation}
and for all but finitely many $\vt\in\vv T$ and all
$\alpha\in\AAA$ there exists a collection $C_{\vt,\alpha}$ of balls $B$ centred in $\Supp$ satisfying the following three conditions:
\begin{equation}\label{e:034}
\Supp\cap \rI_\vt(\alpha,\phi(\vt))\subset \bigcup_{B\in C_{\vt,\alpha}} B,
\end{equation}
\begin{equation}\label{e:035}
\Supp\cap\bigcup_{B\in C_{\vt,\alpha}} B\subset \rI_\vt(\alpha,\phi^+(\vt))
\end{equation}
and
\begin{equation}\label{e:036}
\mu(5B\cap \rI_\vt(\alpha,\phi(\vt)))\le k_\vt \,  \mu(5B).
\end{equation}

\medskip

The following transference theorem is an immediate consequence of \cite[Theorem~5]{Beresnevich-Velani-08-Inhom}.

\begin{theorem}\label{itp}
Suppose that $(\rH,\rI,\Phi)$ satisfies the intersection
property and $\mu$ is contracting with respect to $(\rI,\Phi)$.
Then
\begin{equation*}
\forall\,\phi\in\Phi\quad  \mu(\Lambda_\rH(\phi))=0 \qquad
\Longrightarrow \qquad \forall\,\phi\in\Phi\quad \mu(\Lambda_\rI(\phi)) =0.
\end{equation*}
\end{theorem}

\subsection{Establishing Case B}\label{case_ii}

Recall that our aim is to show that $|\cA_\vf^2(\Psi,\theta)|_m=0$, where $\Psi$ satisfies \eqref{e:001} and (\ref{e:017}). Using \eqref{e:001} and (\ref{e:017}) one readily verifies that
\begin{equation}\label{vb022}
\Psi(\va)<\Psi_0(\va):=\prod_{\stackrel{\scriptstyle i=1 }{a_i\not=0}}^n|a_i|^{-1}
\end{equation}
for all but finitely many $\va\in\ZZ^n$. Therefore,
\begin{equation}\label{vb101}
\cA_\vf^2(\Psi,\theta)\subset \cA_\vf^2(\Psi_0,\theta)
\end{equation} and so
it suffices to show that $|\cA_\vf^2(\Psi_0,\theta)|_m=0$.  With reference to the inhomogeneous transference framework of \S\ref{ITP},
let $\vv T:=(\ZZ_{\ge 0})^n$ and $\AAA:=\ZZ^n\backslash \{0\}\times\Z$. Define the auxiliary function $r:\vv T\to\RR^+$ by setting
\begin{equation}\label{vb021}
r(\vt):=\sqrt{2(n+1)mC_0}\cdot
2^{|\vt|/2}
\end{equation}
where $C_0$ is as in  \eqref{bounds}. Then, given $\varepsilon>0$, $\vt\in\vv T$ and $\alpha=(\vv a,a_0)\in\cA$,   let

\begin{equation}\label{I}
\rI_\vt(\alpha,\varepsilon):=\left\{\vx\in U \,:\,
\left.\begin{array}{rcl}
|a_0+\va\cdot\vf(\vx)+\theta(\vx)|&<&\varepsilon \,\Psi_0(2^\vt)\\[1ex]
|\nabla (\va \cdot\vf(\vx)+\theta(\vx))| &<& \varepsilon \,r(\vt)\\[1ex]
2^{t_i}\le \max\{1,|a_i|\} &<& 2^{t_i+1}\quad (1\le i\le n)\end{array}\right.
\right\}
\end{equation}

\noindent and

\begin{equation}\label{H}
\rH_\vt(\alpha,\varepsilon):=\left\{\vx\in U\,:\,
\left.\begin{array}{rcl}
|a_0+\va\cdot\vf(\vx)|&<&2 \, \varepsilon\, \Psi_0(2^\vt)\\[1ex]
|\nabla (\va\cdot\vf(\vx))| &<& 2\, \varepsilon\, r(\vt)\\[1ex]
|a_i| &<& 2^{t_i+2}\quad(1\le i\le n)
\end{array}\right.\hspace{7.5ex}
\right\}
\end{equation}

\noindent where $2^\vt := (2^{t_1},\ldots,2^{t_n})$. This defines the maps $
\rH $
and $
\rI$ -- see  \eqref{HI}.
Furthermore, given $\delta\in\R$, let $\phi_\delta :\vv T\to\RR^+$ be given by
\begin{equation}\label{vb020}
\phi_\delta(\vt):=2^{\delta |\vt|}    \, ,
\end{equation}
and  let
$$\Phi:=\left\{\phi_\delta \;:\;
0\le\delta<\tfrac{1}{4}\right\}  \, .$$

\noindent For any $\delta\in[0,1/4)$,  it follows  that $$\cA_\vf^2(\Psi_0,\theta)\subset \Lambda_\rI(\phi_\delta)  \,  $$ where
 $\Lambda_\rI(\phi_\delta)$  is the `inhomogenous' $\limsup$ set as defined by (\ref{vb025}).  Therefore, in view of~(\ref{vb101}), to  establish  Case~B it suffices to show that
 \begin{equation}\label{rem6sv}
|\Lambda_\rI(\phi_\delta)|_m=0\quad\text{ for some }\delta\in[0,\tfrac14).
\end{equation}

\noindent With this in mind,  let  $\vv x_0$ be any point in $ U$ at which  $\vv f$ is $l$-non-degenerate  and let $V$ be a sufficiently small open ball centred at $\vv x_0$  such that Corollary~\ref{BKM_corl3_5} and the following statement are valid on $V$.

\begin{theorem}{\bf (\cite[Theorem~1.4]{Bernik-Kleinbock-Margulis-01:MR1829381})} \label{BKM_th1_4}
Let $\vx_0\in U$ and $\vf\,:\,U\to \RR^n$ be $l$-nondegenerate at~$\vx_0$. Then there exists a neighborhood $V\subset U$,
of $\vx_0$ satisfying the following property. For any ball $B\subset V$
there exists $E>0$ such that for any choice of real numbers $\omega,K,T_1,\dots,T_n$ satisfying the inequalities
$$
0<\omega\le 1,\quad T_1,\ldots,T_n\ge 1,\quad
K>0\quad\mbox{ and }\quad\frac{\omega KT_1\cdots T_n}{\max_i T_i}\le 1
$$
the set
$$
S(\omega,K,T_1,\dots,T_n):=\left\{x\in B\,:\, \exists\ \vq\in \ZZ^n \backslash \{0\}\mbox{
such
that}\left. \begin{array}{l} \,\|\vf(\vx)\cdot \vq\|<\omega\\[0.3ex]
|\nabla \vf(\vx)\cdot\vq|<K\\[0.3ex]
\,|q_i|<T_i \quad (1\le i\le n)
\end{array}\right.     \right\}
$$
has $m$-dimensional Lebesgue measure at most $E \, \varepsilon^{\frac{1}{m(2l-1)}}|B|_m$,  where
\begin{equation}\label{vb023}
\varepsilon := \max\left(\omega, \left(\frac{\omega KT_1\cdots T_n}{\max_i
T_i}\right)^{\frac{1}{n+1}}\right).
\end{equation}
\end{theorem}

\noindent Furthermore,  let  $\mu$ be $m$-dimensional Lebesgue measure restricted to $V$. Since $\vv f$ is $l$-non-degenerate almost everywhere, the desired statement \eqref{rem6sv} follows on showing that
\begin{equation}\label{rem6}
\mu(\Lambda_\rI(\phi_\delta))=0\quad\text{ for some }\delta\in[0,\tfrac14)  \, .
\end{equation}
For this, we make use of the Inhomogeneous Transference Principle.  Indeed, suppose for the moment that $(\rH,\rI,\Phi)$ satisfies the intersection
property and $\mu$ is contracting with respect to $(\rI,\Phi)$.  Then, in view of Theorem~\ref{itp}, to establish \eqref{rem6}  it suffices to show that
\begin{equation}\label{vb017}
 \mu(\Lambda_\rH(\phi_\delta))=0   \quad \mbox{for some }\ \delta\in[0,\tfrac14)  \, .
\end{equation}
Armed with Theorem \ref{BKM_th1_4},  it is relatively  painless to establish \eqref{vb017}.  Fix any $\delta\in[0,1/4)$ and notice that in view of  \eqref{H} it follows that
$$
\bigcup_{\alpha\in\cA}\rH_\vt(\alpha,\phi_\delta(\vv t))=S(\omega,K,T_1,\dots,T_n)
$$
with
$$
\omega=2\,\phi_\delta(\vv t)\,  \Psi_0(2^\vt) \, , \quad
K=2\,\phi_\delta(\vv t)\, r(\vt)\quad\text{and}\quad
T_i=2^{t_i+2}\ \  (1\le i\le n).
$$
Using the explicit values of  $\Psi_0(2^\vt)$, $r(\vt)$ and $\phi_\delta(\vv t)$ given by (\ref{vb022}), (\ref{vb021}) and (\ref{vb020}) respectively,  we find that the quantity $\varepsilon$ defined by (\ref{vb023}) satisfies
$$
\varepsilon\ll 2^{-\frac{(1/2-2\delta)}{n+1}|\vt|}.
$$
Therefore, Theorem~\ref{BKM_th1_4} implies that
\begin{equation*}
 \left|  \,  \bigcup_{\alpha\in\cA}  \, \rH_\vt(\alpha,\phi_\delta(\vv t))  \, \right|_m \ll 2^{-\gamma |\vt|}
\end{equation*}
where $\gamma:=\frac{(1/2-2\delta)}{m(n+1)(2l-1)}$ is a positive constant. The  upshot is that
$$
\sum_{\vt\in \vv T}  \,  \left|  \,  \cup_{\alpha\in\cA}  \, \rH_\vt(\alpha,\phi_\delta(\vv t))  \, \right|_m \ll \sum_{\vt\in\ZZ^n}2^{-\gamma|\vt|}<\infty  \, ,
$$
 which together with  the Borel-Cantelli lemma implies the desired  zero  measure  statement
$$
\mu (\Lambda_\rH(\phi_\delta))  = 0  \, .
$$

It remains to verify the intersection and contracting properties.

\subsubsection{Verifying the intersection property}\label{IP}

Let $\vt  \in \vv T$ with  $|\vt| \geq 2  $  and suppose that
$$
\vx\in \rI_\vt(\alpha,\phi_\delta(\vt))\cap \rI_\vt(\alpha',
\phi_\delta(\vt))
$$

\noindent for some distinct $\alpha=(\vv a,a_0)$ and $\alpha'=(\vv a',a_0')$ in $\cA$. Then, by (\ref{I}) and (\ref{H}) we have that

$$
\left\{\begin{array}{rcl} |a_0+\va\cdot \vf (\vx)+\theta(\vx)|&<&\phi_\delta(\vt)\; \Psi_0(2^\vt)\\[0.5ex]
|a_0'+\va'\cdot\vf
(\vx)+\theta(\vx)|&<&\phi_\delta(\vt)\; \Psi_0(2^\vt)
\end{array}\right.
$$
$$
\left\{\begin{array}{rcl} |\nabla(\va\cdot \vf (\vx)+\theta(\vx))|&<&\phi_\delta(\vt)\;  r(\vt)\\[0.5ex]
|\nabla(\va'\cdot\vf
(\vx)+\theta(\vx))|&<&\phi_\delta(\vt)\;  r(\vt)
\end{array}\right.
$$
and
$$
\left\{\begin{array}{rcl} |a_i|&<&2^{t_i+1}\qquad (1\le i\le n)\\[0.5ex]
|a'_i|&<&2^{t_i+1}\qquad (1\le i\le n) \, ,
\end{array}\right.
$$

\noindent where $(a_1,\dots,a_n)=\vv a$ and $(a'_1,\dots,a'_n)=\vv a'$.
Subtracting the first inequality from the second  within each of the above three systems gives

\begin{equation}\label{vb-1more}
\left\{\begin{array}{rcl} |a_0''+\vv a''\cdot\vf (\vx)|&<&2\phi_\delta(\vt)\; \Psi_0(2^\vt)\\[1ex]
|\nabla(\vv a''\cdot \vf (\vx))|&<&2\phi_\delta(\vt)\;  r(\vt)\\[1ex]
|a''_i|&<&2^{t_i+2}\qquad (1\le i\le n)  \; ,
\end{array}\right.
\end{equation}

\noindent where $\vv a''=(a''_1,\dots,a''_n):=\va'-\va$ and $a''_0:=a'_0-a_0$. Regarding the first of the above  inequalities,  by (\ref{vb022}) and the definition of $\Phi$, we have that  $\phi_\delta(\vt)\; \Psi_0(2^\vt)< 2^{-\frac34|\vt|}$. Suppose for the moment that $\vv a''=0$.   Since $\alpha, \alpha' \in\cA $ are distinct, we must have that $a'_0\not=a_0$  and so
$$
|a_0''+\vv a''\cdot\vf (\vx)| =  |a''_0|\ge1  \, .
$$
However,  for any  $\vt $ with  $|\vt| \geq 2  $,  this contradicts the first inequality of (\ref{vb-1more}).
Hence  $\vv a''\not=0$ and  it follows that  $\alpha''\in\cA$. The upshot is  that $\vv x\in \rH_\vt(\alpha'',\phi_\delta(\vv t))$ and therefore  (\ref{e:033}) is satisfied with $\phi^*=\phi_\delta$. This verifies the intersection property.

\subsubsection{Verifying the contracting property}\label{CP}
To start with recall that  $V$ is  a sufficiently small open ball   such that Corollary~\ref{BKM_corl3_5}  is valid  on $5V$.
Thus, there exist positive numbers $H_0 $ and $ C$ such that for any $\vt\in\vv T$ and $\alpha=(\vv a,a_0)\in\AAA$ satisfying $|\va|\ge H_0$ both $a_0+\va\cdot\vf +\theta $ and $|\nabla(\va\cdot \vf +\theta |$ are $(C,\tfrac{1}{ml})$-good on $5V$. In turn, by Lemma~\ref{BKM_lem3_1},  for any $\vt\in\vv T$ and  $\alpha=(\vv a,a_0)\in\AAA$ satisfying $|\va|\ge H_0$ we have that
\begin{equation}\label{vb102}
\text{$\mathbf{F}_{\vt,\alpha}$ is $(C,\tfrac{1}{ml})$-good on $5V$, }
\end{equation}
where $\mathbf{F}_{\vt,\alpha} :U\to\R$  is the function  given by
$$
\mathbf{F}_{\vt,\alpha}(\vx):=\max\Big\{\Psi_0^{-1}(2^\vt)r(\vt)|a_0+\va\cdot\vf(\vx)+\theta(\vx)|,\,
|\nabla(\va\cdot \vf(\vx)+\theta(\vx))|\Big\}.
$$
Notice that the first two inequalities of (\ref{I}) are equivalent to the single inequality $$\mathbf{F}_{\vt,\alpha}(\vx)< \varepsilon \,r(\vt) \, . $$ Therefore, by definition
\begin{equation}\label{vb015}
\rI_\vt(\alpha,\varepsilon)=\left\{\vx\in U \,:\,
\mathbf{F}_{\vt,\alpha}(\vx)< \varepsilon \,r(\vt)
\right\}
\end{equation}
if
\begin{equation}\label{vb015+}
2^{t_i}\le \max\{1,|a_i|\} < 2^{t_i+1}\ \qquad (1\le i\le n).
\end{equation}
Obviously, if (\ref{vb015+}) is not fulfilled then $\rI_\vt(\alpha,\varepsilon)=\emptyset$ irrespective of $\varepsilon$.

\noindent Next,  given $\phi_\delta\in\Phi$ let
$$
\phi^+_\delta:= \phi_{\frac12(\delta+\frac14)}.
$$
Clearly, $\phi^+_\delta$ also lies in $\Phi$. It is easily seen that $\phi_\delta(\vv t)\le\phi^+_\delta(\vv t)$ for all $\vv t\in\vv T$ and therefore
\begin{equation}\label{vb01}
    \rI_\vt(\alpha,\phi_\delta(\vt))\subset \rI_\vt(\alpha,\phi^+_\delta(\vt)).
\end{equation}

We now construct the collection $C_{\vt,\alpha}$ of balls centred in $V$ that satisfy the conditions  (\ref{e:034})--(\ref{e:036}) for an appropriate sequence $k_{\vv t}$. If $\rI_\vt(\alpha,\phi_\delta(\vt))=\emptyset$,  the collection $C_{\vt,\alpha}=\emptyset$ obviously suffices. Thus, we can assume that (\ref{vb015+}) is satisfied and so $\rI_\vt(\alpha,\varepsilon)$ is defined by (\ref{vb015}).
By (\ref{vb022}) and the definition of $\Phi$, it follows that
$$
\rI_\vt(\alpha,\phi^+_\delta(\vt))\subset\{\vv x\in U:|a_0+\va\cdot\vf(\vx)+\theta(\vx)|< 2^{-\frac34|\vt|}\}.
$$
As already pointed out above, $a_0+\va\cdot\vf+\theta$ is $(C,\tfrac{1}{ml})$-good on $5V$ for all  sufficiently large $|\vv a|$.  Therefore, by the  definition of $(C,\alpha)$-good (Definition \ref{def2}) and
Corollary~\ref{BKM_corl3_5+} we have that
\begin{eqnarray*}
|\rI_\vt(\alpha,\phi^+_\delta(\vt))\cap V|_m  & \le  &  |\{\vv x\in V:|a_0+\va\cdot\vf(\vx)+\theta(\vx)| < 2^{-\frac34|\vt|}\}|_m \\[1ex] & \ll  & 2^{-\frac{3|\vv t|}{4ml}}|V|_m  \, ,
\end{eqnarray*}
whenever $|\vv t|$ is sufficiently large. Hence,
\begin{equation}\label{vb-200}
\text{$\rI_\vt(\alpha,\phi^+_\delta(\vt))\not\subset V$  \quad for all  sufficiently large $|\vv t|$.}
\end{equation}

\noindent By (\ref{vb01}) and the fact that $\rI_\vt(\alpha,\phi^+_\delta(\vt))$ is open, for every $\vx\in \Supp\cap \rI_\vt(\alpha,\phi_\delta(\vv t))$ there is a ball $B'(\vv x)$ centred at $\vv x$ such that
\begin{equation}\label{vb011+}
B'(\vv x)\subset \rI_\vt(\alpha,\phi^+_\delta(\vt)).
\end{equation}
On combining  (\ref{vb-200}), (\ref{vb011+}) and the fact that $V$ is bounded, we find that there exists a scaling factor $\tau\ge1$ such that
the ball $B=B(\vv x):=\tau B'(\vv x)$ satisfies
\begin{equation}\label{vb011}
\Supp\cap B(\vv x)\subset \rI_\vt(\alpha,\phi^+_\delta(\vt))\not\supset \Supp\cap 5B(\vv x)
\end{equation}
and
\begin{equation}\label{vb011+2}
5B(\vv x)\subset 5V.
\end{equation}
We now let  $$C_{\vt,\alpha}:=\{B(\vv x)\;:\; \vv x\in \Supp\cap \rI_\vt(\alpha,\phi_\delta(\vt))\}  \, .$$ Then, by construction  and the l.h.s. of (\ref{vb011}), conditions (\ref{e:034}) and (\ref{e:035}) are automatically satisfied. Regarding condition (\ref{e:036}), consider any ball $B\in C_{\vt,\alpha}$. By (\ref{vb015}) and the r.h.s. of~(\ref{vb011}), we have that
\begin{equation}\label{vb013}
\sup_{\vx\in 5B}\mathbf{F}_{\vt,\alpha}(\vx) \, \ge
\, \sup_{\vx\in 5B\cap\Supp} \!\! \mathbf{F}_{\vt,\alpha}(\vx)  \; \ge \; \phi^+_\delta(\vt)  \;  r(\vt).
\end{equation}
On the other hand,
\begin{equation}\label{vb013+}
\sup_{\vx\in 5B\cap
\rI_\vt(\alpha,\phi_\delta(\vt))} \!\!\!\!\!\!\!\!\!\!  \mathbf{F}_{\vt,\alpha}(\vx)\; \le \;  \phi_\delta(\vv t)\; r(\vv t).
\end{equation}
Then, in view of the definitions of $\phi_\delta$, $\phi^+_\delta$ and $r(\vv t)$,  we obtain via  (\ref{vb013}) and (\ref{vb013+}) that
\begin{equation}\label{vb-300}
\sup_{\vx\in 5B\cap
\rI_\vt(\alpha,\phi_\delta(\vt))} \!\!\!\!\!\!\!\!\!\! \mathbf{F}_{\vt,\alpha}(\vx) \ \le  \
2^{-\frac12(\frac1 4-\delta)|\vt|}\sup_{\vx\in 5B}\mathbf{F}_{\vt,\alpha}(\vx).
\end{equation}
Now notice that since \eqref{vb015+} holds,   we have that $|\va|>H_0$ for all $\vv t \in \vv T$ with $|\vv t|$ sufficiently large. Thus, whenever  $|\vv t|$  is sufficiently large,   (\ref{vb102}) is valid which together with  (\ref{vb011+2}) and~(\ref{vb-300}) implies  that
\begin{equation}\label{vb016-}
\begin{array}[b]{rcl}
\Big|5B\cap \rI_\vt(\alpha,\phi_\delta(\vt))\Big|_m
&\le& \Big|\big\{\vv x\in 5B:|\vv F_{\vt,\alpha}(\vv x)|  \le 2^{-\frac12(\frac1 4-\delta)|\vt|}\sup_{\vx\in 5B}\mathbf{F}_{\vt,\alpha}(\vx)\big\}\Big|_m\\[3ex]
&\le &C2^{-\delta^*|\vt|}|5B|_m
\end{array}
\end{equation}
where $\delta^*:=\frac12(\frac14-\delta)\;\tfrac{1}{lm}>0$.
On using the fact that $B$ is centred in $V \subset \Supp$,
 we have that $|5B|_m\le c_m\mu(5B)$ for some constant $c_m$ depending on $m$ only. Hence  (\ref{vb016-}) implies that for all but finitely many $\vt\in\vv T$
$$
\mu(5B\cap \rI_\vt(\alpha,\phi_\delta(\vt))) \, \le  \,
|5B\cap \rI_\vt(\alpha,\phi_\delta(\vt))|_m \, \le \,
c_mC2^{-\delta^*|\vt|}\mu(5B).
$$
This verifies (\ref{e:036}) with $$k_{\vv t}:=c_mC2^{-\delta^*|\vt|}  \, .$$ Furthermore, it is easily seen that the convergence condition \eqref{k_t} is fulfilled. The upshot is that all the conditions of the contracting property are satisfied for the collection  $C_{\vt,\alpha}$  as defined above.

\section{The divergence theory}\label{divergence}

The goal is to prove  Theorems \ref{T3} \&\ \ref{tnice}.  Thus,
throughout  $s > m-1$   and  $\Psi$ is a multivariable approximating function satisfying property $\vv P$ and the divergent sum condition
\begin{equation}\label{e:005+}
 \sum_{\vv a\in\Z^n\nz} |\vv a|  \left(\frac{\Psi(\vv a)}{|\vv a| }\right)^{s+1-m}=\infty.
\end{equation}
Without loss of generality,  we will assume that the vector $\vv v=(v_1, \ldots, v_n)$ appearing in the definition of property $\vv P $ satisfies
\begin{equation}\label{sv2}
v_1=|\vv v|=\max_{1\le i\le n}|v_i|  \, .
\end{equation}

\subsection{Theorem \ref{tnice} \ $\Longrightarrow$ \ Theorem \ref{T3}}\label{aux}

We will need the following technical lemma.

\begin{lemma}\label{measure1}
Let $\mu$ be a finite doubling Borel regular measure on a metric space $(X,d)$ such that $X$ can be covered by a countable collection of arbitrarily small balls.
Let $f:X\to\R^+$ be a uniformly continuous bounded function and let $\nu$ be a measure on $X$ given by
\begin{equation}\label{vbw0}
\nu(A):=\int_A f(x)d\mu(x)
\end{equation}
for every measurable set $A\subset X$.
Let $\{S_Q\}_{Q\in\N}$ be a sequence of measurable subsets of $X$ and $0<\omega<1$ be a constant. Suppose that for every sufficiently small closed ball $B\subset X$
\begin{equation}\label{vbp}
\limsup_{Q\to\infty} \ \mu(S_Q\cap B)\ \le \ \omega\,\mu(B)\,.
\end{equation}
Then for every measurable set $W\subset X$
\begin{equation}\label{vbw9}
\limsup_{Q\to\infty} \ \nu(S_Q\cap W)\ \le \ \omega\,\nu(W)\,.
\end{equation}
\end{lemma}

\proof{}
Let $W$ be any measurable subset of $X$. For every $\ve>0$ and $\delta>0$ there is a finite collection $\cC_{\ve,\delta}$ of disjoint closed balls with radii $<\delta$ such that
\begin{equation}\label{vbw1}
\mu(W\triangle W_{\ve,\delta})<\ve,
\end{equation}
where $E\triangle F:=(E\setminus F)\cup(F\setminus E)$ and $W_{\ve,\delta}:=\bigcup_{B\in\cC_{\ve,\delta}}B$. This is a consequence of \cite[Theorem~2.2.2]{Federer-69:MR0257325} and the discussion of \cite[p.28]{Beresnevich-Dickinson-Velani-06:MR2184760}.
Since $f$ is bounded, there is a constant $C>0$ such that $\nu(A)\le C\mu(A)$ for every measurable set $A$. Therefore, (\ref{vbw1}) implies that
\begin{equation}\label{vbw1+}
\nu(W\triangle W_{\ve,\delta})<C\ve.
\end{equation}

For every $B\in \cC_{\ve,\delta}$ \ let $s_B:=\sup_{x\in B}f(x)$. Since $f$ is bounded, the quantity $s_B$ is finite.
Next, since $f$ is uniformly continuous, for every $\ve>0$ there is a $\delta>0$ such that for every $B\in\cC_{\ve,\delta}$ we have that
\begin{equation}\label{vbw2}
0\le s_B-f(x)<\ve\qquad\text{for all }x\in B\,.
\end{equation}
Since $\cC_{\ve,\delta}$ is finite, property (\ref{vbp}) implies that there is a sufficiently large $Q_0$ such that for all $Q\ge Q_0$ and any $B\in\cC_{\ve,\delta}$ we have that
\begin{equation}\label{vbw3}
\mu(S_Q\cap B)\ \le \ (\omega+\ve)\,\mu(B)\,.
\end{equation}
Then, for $Q\ge Q_0$ it follows that
$$
\begin{array}{rcl}
\nu(S_Q\cap W) & \stackrel{\eqref{vbw1+}}{\le} & C\ve+\sum_{B\in\cC_{\ve,\delta}}\nu(S_Q\cap B)\\[2ex]
& \stackrel{\eqref{vbw0}}{=} & C\ve+\sum_{B\in\cC_{\ve,\delta}}\int_{S_Q\cap B}f(x)d\mu(x)\\[2ex]
& \stackrel{\eqref{vbw2}}{\le} & C\ve+\sum_{B\in\cC_{\ve,\delta}}s_B\int_{S_Q\cap B}d\mu(x)\\[2ex]
& \stackrel{\phantom{\eqref{vbw2}}}{=} & C\ve+\sum_{B\in\cC_{\ve,\delta}}s_B\,\mu(S_Q\cap B)\\[2ex]
& \stackrel{\eqref{vbw3}}{\le} & C\ve+(\omega+\ve)\sum_{B\in\cC_{\ve,\delta}}s_B\mu(B)\\[2ex]
& \stackrel{\eqref{vbw2}}{\le} & C\ve+(\omega+\ve)\sum_{B\in\cC_{\ve,\delta}}\int_B(f(x)+\ve) d\mu(x)\\[2ex]
& \stackrel{\phantom{\eqref{vbw2}}}{=} & C\ve+(\omega+\ve)\int_{W_{\delta,\ve}}(f(x)+\ve) d\mu(x)\\[2ex]
& \stackrel{\phantom{\eqref{vbw2}}}{=} & C\ve+(\omega+\ve)\big(\nu(W_{\delta,\ve})+\ve\mu(W_{\delta,\ve})\big)\\[2ex]
& \stackrel{\eqref{vbw1}\&\eqref{vbw1+}}{\le} & C\ve+(\omega+\ve)\big(\nu(W)+C\ve+\ve(\mu(W)+\ve)\big)\,.
\end{array}
$$
The latter expression tends to $\omega\nu(W)$ as $\ve\to0$. Since $\nu(S_Q\cap W)$ is independent of $\ve$, we obtain (\ref{vbw9}) as required.
\hfill  $\Box$

\bigskip

Let $\vv f:U\to\R^n$ be a map defined on an open set $U\subset\R^m$. Given an $n$-tuple $\vv v=(v_1,\dots,v_n)$ of positive numbers satisfying $v_1+\dots+v_n=n$, $\delta>0$ and $Q>1$, let
$$
\Phi_{\vv v}^{\vf}(Q,\delta)\, :=  \, \Big\{\vx\in U: \exists\ \vv a\in\Z^n\setminus\{0\} \text{ such that }\|\va\cdot\vv f(\vx)\|<\delta Q^{-n}\text{ and } |\vv a|_{\vv v}\le Q\Big\}.
$$

\begin{definition}\label{vnice}\rm
We will say that $\vv f$ is \emph{$\vv v$-nice at $\vv x_0\in U$}\/ if there is a neighborhood $U_0\subset U$  of $\vx_0$ and constants $0<\delta,\omega<1$ such that for any sufficiently small ball $B\subset U_0$ we have that
$$
\limsup_{Q\to\infty}\,|\Phi_{\vv v}^{\vf}(Q,\delta)\cap
B|_m\le \omega|B|_m\,.
$$
The map $\vf$ is said to be \emph{$\vv v$-nice}\/ if it is $\vv v$-nice at almost every point in $U$. Furthermore, $\vf$ is said to be \emph{nice}\/ if it is $\vv v$-nice for all choices of $\vv v$.
\end{definition}

Let $A$ be any Lebesgue measurable subset of $U$. Consider the  measure $ \nu$ given by
$$
\nu(A):=\int_{A} \det G(\vv x)^{1/2}d\vv x,
$$
where $G(\vv x):=\big(g_{i,j}(x)\big)_{1\le i,j\le m}$ with $\vv g_{i,j}:=\partial\vv f/\partial x_i\cdot\partial\vv f/\partial x_j$.   It is well known that the induced measure of a set $S$ on the manifold $\cM$ parameterised by $\vv f$   is given by $\nu(A)$ with $A=\vv f^{-1}(S)$. It is easily verified that
$$
|A|_m=\int_{A}\det G(\vv x)^{-1/2}d\nu(\vv x).
$$
Since $\vv f$ is a Monge  parameterisation, $\det G(\vv x)$ is bounded away from both zero and infinity on a sufficiently small neighborhood of any point $\vv x$. Hence, together with  Lemma~\ref{measure1} we deduce the following statement.

\begin{prop}\label{nice}
Let $\vf:U\to\R^n$ be a $C^2$ parameterisation of a $C^2$ manifold $\cM\subset \R^n$. Let $\vx_0\in U$ and $\vy_0=\vf(\vx_0)$. Then
$\vf$ is $\vv v$-nice at $\vx_0$ if and only if $\cM$ is $\vv v$-nice at $\vy_0$.
\end{prop}

In turn this proposition together with the following lemma implies that non-degenerate manifolds are nice and so Theorem~\ref{T3} is a consequence of Theorem~\ref{tnice}.

\begin{lemma}\label{lem5}
Let $\vv f$ be non-degenerate at $\vv x_0\in U$. Then there is a
ball $B_0\subset U$ centred at $\vx_0$ and a constant $C>0$ such
that for any ball $B\subset B_0$ we have $|\Phi_{\vv v}^{\vv f}(Q,\delta)\cap
B|_m\le C \delta|B|_m$ for all sufficiently large $Q$.
\end{lemma}

\noindent In the case $\vv v=(1,\dots,1)$,  the lemma coincides with  Theorem 2.1 in \cite{Beresnevich-Bernik-Kleinbock-Margulis-02:MR1944505}. For arbitrary $\vv v$, on replacing  the supremum norm by the $\vv v$-quasinorm, the arguments in \cite{Beresnevich-Bernik-Kleinbock-Margulis-02:MR1944505} can be naturally adapted to establish Lemma \ref{lem5}.   The details are left to the energetic reader.

\subsection{Ubiquitous systems in $\RR^m$}\label{us}

The proof of Theorem~\ref{tnice} will make use of the ubiquity framework developed in \cite{Beresnevich-Dickinson-Velani-06:MR2184760}. The framework introduced below is a much simplified version of that in \cite{Beresnevich-Dickinson-Velani-06:MR2184760} and takes into consideration the specific application that we have in mind.

Throughout, balls in $\R^m$ are assumed to be defined in terms of the supremum norm $|\cdot|$.
Let $U$ be a ball in $\RR^m$ and $\RRR=(R_\alpha)_{\alpha\in J}$ be a family of subsets
$R_\alpha\subset \R^m$ indexed by a countable set $J$. The sets $\ra$ are referred to as \emph{resonant sets}. Throughout, $\rho\;:\;\RR^+\to\RR^+$
will denote a function such that
$\rho(r)\to0$ as $r\to\infty$.  Given a set $A\subset U$, let
$$
\Delta(A,r):=\{\vx\in U\;:\; \dist(\vx,A)<r\}
$$
where $\dist(\vv x,A):=\inf\{|\vv x-\vv a|: \vv a\in A\}$. Next, let
$\beta\;:\;J\to \RR^+\;:\;\alpha\mapsto\beta_\alpha$ be a positive
function on $J$. Thus the function $\beta$ attaches a `weight'
$\beta_\alpha$ to the set $R_\alpha$. We will assume that for every
$t\in \NN$ the set $J_t=\{\alpha\in J: \beta_\alpha\le 2^t\}$ is
finite.

\bigskip

\noindent\textbf{The intersection conditions:} There exists a constant $\gamma$ with $ 0 \leq \gamma \leq m$ such that for any sufficiently large $t$ and for any
$\alpha\in J_t$, $c\in\ra$ and $0< \lambda \le \rho(2^t)$ the
following conditions are satisfied:
\begin{equation}\label{i1}
\big|B(c, {\mbox{\small
$\frac{1}{2}$}}\rho(2^t))\cap\Delta(\ra,\lambda)\big|_m \geq c_1 \,
|B(c,\lambda)|_m\left(\frac{\rho(2^t)}{\lambda}\right)^{\gamma}
\end{equation}
\begin{equation}\label{i2}
\big|B\cap
B(c,3\rho(2^t))\cap\Delta(\ra,3\lambda)\big|_m\leq c_2 \,
|B(c,\lambda)|_m\left(\frac{r(B)}{\lambda}\right)^{\gamma} \
\end{equation}
where $B$ is an arbitrary ball centred on a resonant set with radius
$r(B)\le 3 \, \rho(2^t)$. The   constants $c_1$ and $ c_2$  are
positive and absolute. The constant $\gamma$ is referred to as the \emph{common dimension} of $\RRR$.

\begin{definition}\rm
Suppose that there exists a ubiquitous function $\rho$ and an
absolute constant $k>0$ such that for any ball $B\subseteq U$
\begin{equation}\label{coveringproperty}
\liminf_{t\to\infty} \left|\bigcup_{\alpha\in
J_t}\Delta(R_\alpha,\rho(2^t))\cap B\right|_m\ \ge \ k\,|B|_m\,.
\end{equation}
Furthermore, suppose that the intersection conditions \eqref{i1} and \ref{i2}  are satisfied. Then
the system $(\RRR, \beta)$ is called \emph{locally ubiquitous in $U$
relative to $\rho$.}
\end{definition}

Let $(\RRR,\beta)$ be a ubiquitous system in $U$ relative to $\rho$ and
$\phi$ be an approximating function. Let $\Lambda(\phi)$ be the set
of points $\vv x\in U$ such that the inequality
\begin{equation}\label{vb+}
 \dist(\vv x,R_{\alpha})<\phi(\beta_\alpha)
\end{equation}
holds for infinitely many $\alpha\in J$.

\begin{lemma}[Ubiquity Lemma]\label{lem1_intro}
Let $\phi$ be an approximating function and $(\RRR, \beta)$ be a
locally ubiquitous system in $U$ relative to $\rho$. Suppose that there
is a $\lambda\in\RR$, $0<\lambda<1$ such that
$\rho(2^{t+1})<\lambda\rho(2^t)$ for all $t\in\N$. Then for any
$s>\gamma$
\begin{equation}\label{uld}
\HHH^s(\Lambda(\phi))=\HHH^s(U)\qquad\text{if}\qquad
\sum_{t=1}^\infty\frac{\phi(2^t)^{s-\gamma}}{\rho(2^t)^{m-\gamma}}=\infty.
\end{equation}
\end{lemma}

\vspace*{2ex}

\noindent\emph{Remark.}  When $s>m$, we have that $\HHH^s(\Lambda(\phi))=\HHH^s(U)=0$ and the lemma is trivial.  In the case $s=m$ it is a consequence of \cite[Corollary~2]{Beresnevich-Dickinson-Velani-06:MR2184760} and in the case $s<m$ it is a consequence of
\cite[Corollary~4 ]{Beresnevich-Dickinson-Velani-06:MR2184760}.

\subsection{The appropriate ubiquitous system for Theorem \ref{tnice}  }

Recall that $\vv f=(f_1,\dots,f_n):U\to\R^n$ is a $\vv v$-nice $C^2$ map satisfying (\ref{bounds}), where $U$ is a ball in $\R^m$. Also recall that $\theta:U\to\R$ is a $C^{(2)}$ function. Let $\FFF_n$ denote the set of all functions $F:U\to\R$ given by
$$
F(\vv x)=a_0+a_1f_1(\vx)+a_2f_2(\vx)+\ldots+a_nf_n(\vx)\,,
$$
where $a_0,\ldots,a_n$ are integer coefficients not all zero.
Given $F\in\FFF_n$, let
\begin{equation}\label{sv3}
\tilde R_F:=\{\vv x\in U: F(\vv x)+\theta(\vv x)=0\}\qquad\text{and}\qquad
H_{\vv v}(F):=\max_{1\le i\le n}|a_i|^{1/v_i}.
\end{equation}

\noindent The key to establishing Theorem~\ref{tnice} is the following ubiquity statement.   With reference to the abstract  setup of  \S\ref{us},   the indexing set $J=\FFF_n$ and so $F$ plays the role of $\alpha\in J$.

\begin{prop}\label{prop1}
Let $\vv x_0\in U$ be such that $\vv f$ is $\vv v$-nice at $\vv x_0$. Then there is a neighborhood $U_0$ of $\vv x_0$, constants $\kappa_0>0$ and $\kappa_1>1$ and a collection $\RRR:=\big(R_F\big)_{F\in \FFF_n}$ of sets $R_F\subset \tilde R_F\cap U_0$ such that the system $(\RRR,\beta)$, where
$$
\beta\, : \, \FFF_n\to\R^+\, :\, F\mapsto\beta_F:=\kappa_0H_{\vv v}(F),
$$
is locally ubiquitous in $U_0$ relative to $\rho(r):=\kappa_1r^{-n-v_1}$ with common dimension $\gamma:=m-1$.
\end{prop}

The sets $\tilde R_F$ are essentially the appropriate resonant sets. However,
to ensure that the intersection conditions associated with ubiquity are satisfied,
in particular, the lower bound condition (\ref{i1}),
we cannot in general work with the sets $\tilde R_F$ directly\footnote{In various previous applications of ubiquity  to approximation
problems on manifolds the intersection conditions have not always been explicitly addressed.  Indeed, it is not clear in some instances whether or not  the  authors have defined $\tilde R_F$  to be the resonant sets.}.
To illustrate this, consider the following explicit  examples.

\bigskip

\noindent\emph{Example 1.} Let $m=2$, $n=3$, $U=\{(x_1,x_2)\in\R^2:x_1^2+x_2^2<1\}$ and $f(x_1,x_2)=\sqrt{1-x_1^2-x_2^2}$.
It is easily seen that for most choices of $F$ the intersection conditions are satisfied with $\gamma=1$. However,
when $\vv a=(-1,0,0,1)$ and so $F=f-1$, we have that $\tilde R_F=\{(0,0)\}$. Then the l.h.s. of (\ref{i1}) is comparable to $\lambda^2$, while the r.h.s. of (\ref{i1}) is comparable to $\lambda\rho(2^t)$.  Thus (\ref{i1}) is violated.

\bigskip

\noindent\emph{Example 2.} Let $m=2$, $n=3$, $U=(\alpha,\alpha+1)^2$ with $\alpha$ a Liouville number and $f(x_1,x_2)=x_1^2+x_2^2$.
As in the above example it is easily seen that for most choices of $F$ the intersection conditions are satisfied with $\gamma=1$. Since $\alpha$ is Liouville, for any real   $v$ we have that $|\alpha-p/q|<q^{-v}$ for infinitely many rationals $p/q$  $(q>0)$. Consider $\vv a=(-2p, q,q,0)$ if $\alpha-p/q<0$ and $\vv a=(-2(p+q),q,q,0)$ if $\alpha-p/q>0$. It is a simple matter to verify that $\tilde R_F$ is a  line segment of length comparable to  $|\alpha -p/q|< q^{-v}$.  Then the l.h.s. of (\ref{i1}) is comparable to $\lambda(\lambda+ q^{-v})$, while the r.h.s. of (\ref{i1}) is comparable to $\lambda\rho(2^t)$. For large enough $v$, the upshot is that  (\ref{i1}) is violated.

\bigskip

\noindent The upshot  is that the sets  $\tilde R_F$  need to be modified in an appropriate manner to yield the resonant sets $R_F$ -- namely via the `trimming' procedure described in \S\ref{trs} below\footnote{The trimming procedure can be replicated to address  the oversights alluded in the   previous footnote.}.

\subsubsection{Proof of Theorem~\ref{tnice} modulo Proposition~\ref{prop1}}

Fix $\vv x_0\in U$ such that $\vv f$ is $\vv v$-nice at $\vv x_0$ and let $U_0$ be as in Proposition~\ref{prop1}. Since $\vv f$ is $\vv v$-nice (i.e. $\vv f$ is $\vv v$-nice at almost every point in $U$), it suffices to prove that
\begin{equation}\label{ohyes}
\HHH^s(\cA_{\vv f}(\Psi,\theta)\cap U_0)=\HHH^s(U_0)\,.
\end{equation}

\bigskip

\noindent With reference to \S\ref{us},  let $U = U_0$ and
$$
\phi : r \to \phi(r):= (2 n C_0)^{-1} (\kappa_0^{-1}r)^{-v_1} \, \psi(\kappa_0^{-1}r)  \, .
$$
Here  the approximating function $\psi$ and the  vector   $\vv v=(v_1,\dots,v_n)$ are  associated with the fact that  $\Psi$ is a multivariable  approximation function  satisfiing  property $\vv P$.  Our first goal is to show that
\begin{equation}\label{goal1}
\Lambda(\phi)\subset\cA_{\vv f}(\Psi,\theta).
\end{equation}
Let $\vv x=(x_1,\dots,x_m)\in\Lambda(\phi)$. By definition, $ \Lambda(\phi) $ is a subset of $U_0$ and inequality (\ref{vb+}) is  satisfied for infinitely many $F=a_0+a_1f_1+\dots+a_nf_n\in\FFF_{n}$  -- recall that we have identified  $\alpha$ with $F$ and $J$ with $\FFF_{n}$.   Now fix such a function $F$. Then,   by the definition of $\beta$ and the properties of $R_F$ within Proposition~\ref{prop1}, there exists a point $\vv z=(z_1,\dots,z_m)\in U_0$ such that $F(\vv z)+\theta(\vv z)=0$ and
\begin{equation}\label{sv1}
    |\vv x-\vv z|<\phi(\kappa_0H_{\vv v}(F)).
\end{equation}
Thus,  by  the  Mean Value Theorem   it follows that there exists some $\tilde\vx  \in U_0$ such that
$$
\begin{array}{rcl}
|F(\vx)+\theta(\vx)| & = & \Big|\sum_{i=1}^m
\tfrac{\partial}{\partial x_i}(F+\theta)(\tilde\vx)(x_i- z_{i})\Big|\\[2ex]
& \le & |\vv x-\vv z|\,\sum_{i=1}^m
\Big|\tfrac{\partial}{\partial x_i}\big(\sum_{j=1}^n a_jf_j+\theta\big)(\tilde\vx)\Big|\\[2ex]
& \stackrel{\eqref{bounds}}{\le} & 2 n C_0 \ |\vx-\vv z|\,
\max_{1\le j\le n}|a_j| \\[2ex]
& \stackrel{\eqref{sv1}}{\le} &  2 n C_0  \  \phi(\kappa_0H_{\vv v}(F))\,
\max_{1\le j\le n}|a_j| \\[2ex]
& \stackrel{\eqref{sv2}+\eqref{sv3}}{\le} & 2 n C_0 \ \phi(\kappa_0H_{\vv v}(F))\,
H_{\vv v}(F)^{v_1} \\[2ex]
& \stackrel{}{\le} &  \psi(H_{\vv v}(F))  \ = \ \Psi(\vv a)  \, .
\end{array}
$$
The upshot is that there are infinitely many $F\in\FFF_n$ satisfying the above inequalities. This verifies (\ref{goal1}) and  together with  Lemma~\ref{lem1_intro}  implies  \eqref{ohyes} as long as the
 sum in (\ref{uld}) diverges.  We now verify this divergent condition.  Recall that $\gamma:=m-1$ and so
\begin{equation}\label{uld2}
\sum_{t=1}^\infty\frac{\phi(2^t)^{s-m+1}}{\rho(2^t)}  \ \asymp  \
\sum_{t=1}^\infty\frac{(2^{-v_1t}\psi(\kappa_0^{-1}2^t))^{s-m+1}}{2^{-(n+v_1)t}}  \, .
\end{equation}
On using the fact that $v_1+\dots+v_n=n$, it follows that for any $t\in\N$ the number of $\vv a\in\Z^n$ such that $\kappa_02^{t}<|\vv a|_{\vv v}\le \kappa_02^{t+1}$ is comparable to $ 2^{nt}$. Also, by (\ref{sv2})  we have that $|\vv a|\asymp2^{v_1t}$ whenever $\kappa_02^{t}<|\vv a|_{\vv v}\le \kappa_02^{t+1}$. Therefore,
\begin{equation}\label{uld3}
\text{r.h.s. of (\ref{uld2}) }  \ \asymp  \  \sum_{t=1}^\infty\ \ \sum_{\kappa_02^{t}<|\vv a|_{\vv v}\le \kappa_02^{t+1}}|\vv a|\left(\frac{\psi(\kappa_0^{-1}2^t)}{|\vv a|}\right)^{s-m+1}.
\end{equation}
Next, since  $\psi$ is decreasing, it follows that $\psi(\kappa_0^{-1}2^t)\ge\psi(|\vv a|_{\vv v})=\Psi(\vv a)$ whenever $\kappa_02^{t}<|\vv a|_{\vv v}\le \kappa_02^{t+1}$. Therefore,
\begin{eqnarray*}
\text{r.h.s. of (\ref{uld3}) }  & \gg  &  \sum_{t=1}^\infty \sum_{\kappa_02^{t}<|\vv a|_{\vv v}\le \kappa_02^{t+1}}|\vv a|\left(\frac{\Psi(\vv a)}{|\vv a|}\right)^{s-m+1} \\ \nonumber & \asymp & \sum_{\vv a\in\Z\nz}|\vv a|\left(\frac{\Psi(\vv a)}{|\vv a|}\right)^{s-m+1} \ \ \stackrel{\eqref{e:005+}}= \ \infty.
\end{eqnarray*}
This completes the proof of Theorem~\ref{tnice} modulo  Proposition~\ref{prop1}.

\subsubsection{The resonant sets \label{trs}}

As already mentioned, the sets $\widetilde R_F$ given by (\ref{sv3}) are essentially the appropriate resonant sets. However, to ensure that the intersection conditions  associated with ubiquity are satisfied, these sets require modification.  Essentially, we impose the condition that
\begin{equation}\label{extra}
 |\tfrac{\partial}{\partial x_1} (F+\theta)(\vv x)| \; > \; p \, |\nabla (F+\theta)(\vv x)|\qquad  \text{for all $\vv x\in U_0$}
\end{equation}
for some  fixed  $p\in(0,1)$. In what follows the projection map $\pi:\R^m\to\R^{m-1}$ will be given by
\begin{equation}\label{pi}
\pi(x_1,x_2,\dots,x_m)=(x_2,\dots,x_m)\,.
\end{equation}

\begin{prop}\label{prop3}
Let $\rho$ and $\beta$ be as in Proposition~\ref{prop1}. Let $U_0 $ be any open subset of $U$ and  $p\in(0,1)$. For $F\in\FFF_n$ let
\begin{equation}\label{ressets}
\widetilde V:=\pi(\widetilde R_F\cap U_0),\qquad
V:=\hspace*{-2ex}\bigcup_{3\rho(\beta_F)\text{-balls } B\,\subset\, \widetilde V}\hspace*{-4ex}\tfrac12B
\end{equation}
and
\begin{equation}\label{ressets2}
    R_F:=\left\{\begin{array}{cccl}
  \pi^{-1}(V)\cap\widetilde R_F  &&& \text{ if $F$ satisfies \eqref{extra}}\\[1ex]
  \emptyset &&& \text{ otherwise}
           \end{array}
\right.
\end{equation}
where $3\rho(\beta_F)$-balls are open balls in $\R^{m-1}$ of radius $3\rho(\beta_F)$.
Then, $R_F$ satisfies the intersection conditions \eqref{i1} and \eqref{i2} with
$$
c_1:=2^{-2m+3}v_m^{-1}\qquad\text{ and }\qquad c_2:=3m2^m(p\,v_m)^{-1}\,,
$$
where $v_m$ is the volume of an $m$-dimensional ball of unit radius.
\end{prop}

\bigskip

\proof{}
Let $t\in\N$, $F\in\FFF_n$ and $\beta_F\le 2^t$. In view of (\ref{extra}) the gradient of $F+\theta$ never vanishes on $U_0$ and therefore the set $\widetilde R_F\cap U_0:=\{\vv x\in U_0:F(\vv x)+\theta(\vv x)=0\}$ is a regular $C^{(2)}$ submanifold of $U_0$ of dimension $(m-1)$. This is a well known fact from differential geometry -- see, for example \cite[Theorem~1.13]{Olver-93:MR1240056}. Furthermore,
(\ref{extra}) together with the Implicit Function Theorem implies that $R_F\cap U_0$ can be defined as the graph $G_g(\widetilde V)$ of a $C^{(2)}$ function $g:\widetilde V\to\R$, where
\begin{equation}\label{1more}
G_g(S):=\{(g(x_2,\dots,x_m),x_2,\dots,x_m):(x_2,\dots,x_m)\in S\}
\end{equation}
for $S\subseteq \widetilde V$.
Then, by the definition of $R_F$, we have that $R_F=G_g(V)$.
If $R_F$ happens to be empty, the intersection conditions \eqref{i1} and \eqref{i2} are trivially satisfied. Otherwise, $R_F\not=\emptyset$ and we proceed as follows.

Given $r>0$ and a set $A\subset\R^{m}$, let
$$
\Delta_1(A,r):=\{\lambda\vv e_1+\vv x:|\lambda|\le r,\ \vv x\in A\},
$$
where $\vv e_1:=(1,0,\dots,0)\in\R^m$. By the definition of $g$,
$$
(F+\theta)(g(x_2,\dots,x_m),x_2,\dots,x_m)=0\qquad\text{for all }\quad (x_2,\dots,x_m)\in\tilde V\,.
$$
Then differentiating this identity and using (\ref{extra}), we obtain that
\begin{equation}\label{extra2}
    |\nabla g(x_2,\dots,x_m)|\le p^{-1}\quad\text{for all } \ (x_2,\dots,x_m)\in \widetilde V.
\end{equation}
We now show  that
\begin{equation}\label{vbx1}
\Delta_1(R_F,\eta) \, \subset \, \Delta(R_F,\eta)  \, \subset  \,  \Delta_1(\widetilde R_F\cap  U_0,\eta m p^{-1})\quad\text{for any $\eta\le 3\rho(\beta_F)$}.
\end{equation}
Indeed, the l.h.s. of (\ref{vbx1}) is a straightforward consequence of the definitions of $\Delta(A,r)$ and $\Delta_1(A,r)$. To prove the r.h.s. of (\ref{vbx1}) take any $\vv z\in \Delta(R_F,\eta)$. Then there exists $\vv x\in R_F$ such that $\dist(\vv z,\vv x)<\eta$. By the definition of $R_F$ and $V$, we have that $\pi\vv x\in \tfrac12B$ for some $3\rho(\beta_F)$-ball $B\subset \widetilde V$. Hence, $B(\pi\vv x,3\rho(\beta_F))\subset B\subset \widetilde V$. Since $\dist(\pi\vv z,\pi\vv x)\le\dist(\vv z,\vv x)<\eta\le 3\rho(\beta_F)$, we have that $\pi\vv z\in \widetilde V$. Then, on making use of the Triangle Inequality and the Mean Value Theorem we find that
$$
|z_1-g(\pi\vv z)|\stackrel{\eqref{1more}}{=} |z_1-x_1+g(\pi\vv x)-g(\pi\vv z)|\le \eta+|g(\pi\vv x)-g(\pi\vv z)|\stackrel{\eqref{extra2}}{\le} \eta mp^{-1} \, .
$$
This verifies the r.h.s. of (\ref{vbx1}).  We are now in the position to establish  the intersection conditions \eqref{i1} and \eqref{i2}.

\medskip

\noindent\emph{The lower bound condition}\/. Let $\vv c\in R_F$ and $0<\lambda\le\rho(2^t)$. Since $\rho$ is decreasing, we have that $\rho(2^t)\le\rho(\beta_F)$.
Then, by (\ref{vbx1}), we find that
\begin{equation}\label{vbx3}
B(\vv c,\tfrac12\rho(2^t))\cap\Delta(R_F,\lambda) \ \supset \ B(\vv c,\tfrac12\rho(2^t))\cap\Delta_1(R_F,\lambda)
 \ \supset \ \Delta_1(G_g(W),\lambda ),
\end{equation}
where $W:=\pi(B(\vv c,\tfrac12\rho(2^t)))\cap V$. Since $\vv c\in R_F$, we have that $\pi\vv c\in V$ and therefore there exists a $3\rho(\beta_F)$-ball $B\subset \widetilde V$ such that $\pi\vv c\in \tfrac12B$. Hence, since $3\rho(\beta_F)\ge\rho(2^t)$ and $\pi\vv c\in \tfrac12B\subset V$, the set $\pi(B(\vv c,\tfrac12\rho(2^t)))\cap \tfrac12B$ contains a ball of radius $\tfrac14\rho(2^t)$  and therefore
$$
|\pi(B(\vv c,\tfrac12\rho(2^t)))\cap \tfrac12 B|_{m-1}\ge (\tfrac14\rho(2^t))^{m-1}v_{m-1}\ge (\tfrac14\rho(2^t))^{m-1}.
$$
Consequently, $|W|_{m-1}\ge (\frac14\rho(2^t))^{m-1}$.
Finally using (\ref{vbx3}) and Fubini's theorem gives
$$
\begin{array}{rcl}
\big|B(\vv c,\frac12\rho(2^t))\cap\Delta(R_F,\lambda)\big|_m & \ge  & \big|W\big|_{m-1}\  2\lambda
 \ \ \ge  \ \ (\frac14\rho(2^t))^{m-1}\ 2\lambda \\[3ex]
 & =  & c_1\,|B(\vv c,\lambda)|_{m}\left(\dfrac{\rho(2^t)}{\lambda}\right)^{m-1}  \, .
 \end{array}
$$

\medskip

\noindent\emph{The upper bound condition}\/. Take any $\vv c\in R_F$, any positive $\lambda\le\rho(2^t)$ and any ball $B$ with radius $r(B)\le 3\rho(2^t)$. Since $\rho$ is decreasing, we also have that $\rho(2^t)\le\rho(\beta_F)$.
Then, by~(\ref{vbx1}), we find that
\begin{equation}\label{vbx2}
\begin{array}{rcl}
B\cap B(\vv c,3\rho(2^t))\cap\Delta(R_F,3\lambda) & \subset & B\cap B(\vv c,3\rho(2^t))\cap\Delta_1(\widetilde R_F\cap  U_0,3\lambda m p^{-1})\\[2ex]
 & \subset & \Delta_1(G_g(W'),3\lambda m p^{-1}),
\end{array}
\end{equation}
where $W':=\pi(B\cap B(\vv c,3\rho(2^t))\cap\widetilde R_F\cap  U_0)$. Clearly, $\diam\, W'\le 2r(B)$. Therefore, using~(\ref{vbx2}) and Fubini's theorem gives
$$
\begin{array}{rcl}
\big|B\cap B(\vv c,3\rho(2^t))\cap\Delta(R_F,3\lambda)\big|_m & \le  & \big|W'\big|_{m-1}\ 6\lambda m p^{-1}
\ \ \le  \ \ (2r(B))^{m-1}\ 6\lambda m p^{-1}\\[3ex]
 & =  & c_2\,|B(\vv c,\lambda)|_{m}\left(\dfrac{r(B)}{\lambda}\right)^{m-1}  \, .
 \end{array}
$$
\hfill  $\Box$

\bigskip

\subsection{Proof of Proposition~\ref{prop1}}

Let $\vv x_0\in U$ be such that $\vv f$ is $\vv v$-nice at $\vv x_0$ and let $U_0$ be the neighborhood of $\vv x_0$ that arises from Definition~\ref{vnice}. Without loss of generality, we will assume that $U_0$ is a ball satisfying
\begin{equation}\label{U_0}
\diam\, U_0\le \big(2nm(n+1)C_0\delta^{-n}\big)^{-1}\,,
\end{equation}
where $\delta$ is as in Definition~\ref{vnice} and $C_0$ is as in (\ref{bounds}).
We shall show that there are constants $\kappa_0>0$ and $\kappa_1>1$ and a value for $p$  associated with \eqref{extra} such that the collection $(R_F)_{F\in\FFF_n}$ given by (\ref{ressets2})
satisfies the statement of Proposition~\ref{prop1}. In view of  Proposition~\ref{prop3}, the intersection conditions (\ref{i1}) and (\ref{i2}) are then automatically satisfied. Thus, to establish ubiquity all that remains is to verify the  measure theoretic `covering'   condition~(\ref{coveringproperty}).

Let $B\subset U_0$ be an arbitrary ball and $t$ be a sufficiently large integer. Let $$Q=2^t.$$
  By Definition~\ref{vnice}, for some fixed $\delta,\omega\in(0,1)$ we have that
$$
\limsup_{Q\to\infty}\,|\Phi_{\vv v}^{\vf}(Q,\delta)\cap
\tfrac12B|_m\le \omega|\tfrac12B|_m\,.
$$
Therefore, for sufficiently large $Q$ we have that
$$
|\tfrac12B\setminus\Phi_{\vv v}^{\vv f}(Q,\delta)|_m\ge \tfrac12(1-\omega) |\tfrac12B|_m = 2^{-m-1}(1-\omega) |B|_m.
$$
Therefore, if we can show that
\begin{equation}\label{2more}
\tfrac12B\setminus\Phi_{\vv v}^{\vv f}(Q,\delta)\subset \bigcup_{\stackrel{\scriptstyle F\in\FFF_n}{\beta_F\le Q}}\Delta(R_F,\rho(Q))\cap B
\end{equation}
then (\ref{coveringproperty}) would follow as required. With this in mind,
let $$\vv x\in \tfrac12B\setminus\Phi_{\vv v}^{\vv f}(Q,\delta)$$ and consider the system of inequalities
\begin{equation}\label{e:013}
\left\{\begin{array}{rcl}
|a_nf_n(\vv x)+\ldots+a_1f_1(\vv x)+a_0| & < & Q^{-n}\\[1ex]
|a_i| & \le & Q^{v_i}\qquad (1\le i\le n).
\end{array}\right.
\end{equation}
The set of $(a_0,\dots,a_n)\in\R^{n+1}$ satisfying (\ref{e:013}) gives rise to  a convex body $D$ in $\RR^{n+1}$ which is symmetric about the origin. Let $\tau_0,\ldots,~\tau_{n+1}$ be the successive minima of $D$. By
definition, $\tau_1\le\tau_2\le\ldots\le~\tau_{n+1}$. Since $\vv x\not\in \Phi_{\vv v}^{\vv f}(Q,\delta)$, we have that
$\tau_1\ge\delta$.
By Minkowski's theorem on successive minima \cite{Cassels-1957}, we have that
$$\tau_1\cdots\tau_{n+1}\operatorname{Vol}(D)\le 2^{n+1}\, . $$ In view of the fact that $v_1+\dots+v_n=n$ we find that $\operatorname{Vol}(D)=2^{n+1}$. Therefore, $\tau_1\cdots\tau_{n+1}\le 1$, whence
$$
\tau_{n+1}\le (\tau_1\cdot \tau_2\cdots\tau_n)^{-1}<\delta^{-n}\,.
$$
By the definition of $\tau_{n+1}$, there are linearly
independent vectors $\vv a_j=(a_{j,0},\dots,a_{j,n})\in\Z^{n+1}$ $(0\le
j\le n)$ such that the functions $F_j$ given by
$$F_j(\vv x):=a_{j,n}f_n(\vv x)+\ldots+a_{j,1}f_1(\vv x)+a_{j,0}$$ satisfy
\begin{equation}\label{e:014}
\left\{\begin{array}{lcl}
|F_j(\vv x)|&\le& C_2Q^{-n}\\[1ex]
|a_{j,i}|&\le& C_2Q^{v_i}\qquad (1\le i\le n),
\end{array}\right.
\end{equation}
where
\begin{equation}\label{C_2}
    C_2:=\delta^{-n}.
\end{equation}

\noindent The  next step is to construct a linear combination of $F_j$ which gives rise  to a resonant set $R_F$ with  $\vv x$ lying  within a sufficiently  small neighborhood of $R_F$. With this in mind, consider the following system of linear equations
\begin{equation}\label{e:015}
\left\{\begin{array}{rcl}
\eta_0F_0(\vv x)+\ldots+\eta_{n}F_{n}(\vv x)+\theta(\vv x)&=&0\\[1.5ex]
\eta_0\frac{\partial}{\partial x_1}F_0(\vv x)+\ldots+\eta_{n}\frac{\partial}{\partial x_1}F_{n}(\vv x)+\frac{\partial}{\partial x_1}\theta(\vv x)&=&Q^{v_1}+\sum_{i=0}^n|\frac{\partial}{\partial x_1}F_i(\vv x)|\\[1.5ex]
\eta_0a_{0,j}+\ldots+\eta_na_{n,j}&=&0\qquad (2\le j\le n).
\end{array}\right.
\end{equation}
Using the fact that $f_1(\vv x)=x_1$, it is readily verified that the determinant of this system is equal to
$\det(a_i^{(j)})_{0\le i,j\le n}$. The latter is non-zero since $\vv a_0,\dots,\vv a_n$ are linearly independent. Therefore, the
system~(\ref{e:015}) has a unique solution
$\eta_0,\ldots,\eta_{n}$. For the integers $t_i:=\lfloor\eta_i\rfloor$ we have that
\begin{equation}\label{t}
|t_i-\eta_i|<1\qquad (0\le i\le n).
\end{equation}
Let
$$
  F(\vv x) := t_0F_0(\vv x)+\ldots+t_{n}F_{n}(\vv x)
= a_nf_n(\vv x)+\ldots+a_1f_1(\vv x)+a_0,
$$
where $a_i:=t_0a_{0,i}+\ldots+t_{n}a_{n,i}$. We claim that $F$ satisfies (\ref{extra}), the height condition $\beta_F\le Q$ and moreover $\vv x\in \Delta(R_F,\rho(Q))$. Thus (\ref{2more}) follows and we are done.

\bigskip

\noindent\textbf{Verifying the height condition:}
By making use of (\ref{e:014}), (\ref{e:015}) and (\ref{t}),  we find that
\begin{equation}\label{vb4}
    |a_j|\le (n+1)C_2 Q^{v_i} \qquad(2\le j\le n)
\end{equation}
and
\begin{equation}\label{vb1}
 |F(\vv x)+\theta(\vv x)|\le (n+1)C_2 Q^{-n}.
\end{equation}
Using the second equation of (\ref{e:015}), we find that
\begin{equation}\label{vb2}
\Big|\tfrac{\partial}{\partial
x_1}\big(F+\theta\big)(\vv x)\Big|\ge Q^{v_1}.
\end{equation}
In particular, this means that $F$ is not identically zero and so $F\in\FFF_n$.
Next, using (\ref{bounds}), (\ref{e:014}) and the assumption that $v_1=|\vv v|$ we find that $$\Big|\frac{\partial}{\partial x_1}F_i(\vv x)\Big|\le nC_0Q^{v_1}  \quad\text{ for  all }  i=0,\dots,n. $$ Together with (\ref{e:015}) and (\ref{t}),  this implies that
\begin{equation}\label{vb3}
\left|\tfrac{\partial}{\partial x_1}\big(F+\theta\big)(\vv x)\right|\le(2nC_0+1) Q^{v_1}.
\end{equation}
Furthermore, since $\vv f$ is  a Monge parameterisation we have that
$$
a_1=\tfrac{\partial}{\partial x_1}(F+\theta)(\vv x)-\tfrac{\partial}{\partial x_1}\theta(\vv x)-\sum_{j=2}^na_j\tfrac{\partial}{\partial x_1}f_j(\vv x) \, .
$$
Then, on  using (\ref{bounds}), (\ref{vb4}) and (\ref{vb3}) we obtain that
\begin{equation}\label{a1}
|a_1|\le C_3Q^{v_1},\qquad\text{where}\quad C_3:=(n+3)^2C_0C_2\,.
\end{equation}
This together with (\ref{vb4}) and (\ref{vb2}) gives that
\begin{equation*}
    \kappa_0^*\,Q\le \beta_F:=\kappa_0H_{\vv v}(F)\le Q \, ,
\end{equation*}
for some explicitly computable constant $\kappa_0,\kappa_0^*>0$ depending only on $\vv v$, $n$, $C_0$ and $C_2$.

\bigskip

\noindent\textbf{Verifying condition (\ref{extra}):}
In view of  Taylor's formula,  for any $\vv y\in U_0$ we have that
\begin{equation}\label{new1}
\left|\tfrac{\partial}{\partial x_1}(F+\theta)(\vv y)\right| \ \ge\ \left|\tfrac{\partial}{\partial x_1}(F+\theta)(\vv x)\right|-\sum_{i=1}^m \left|\tfrac{\partial^2}{\partial x_1\partial x_i}(F+\theta)(\tilde {\vv y})(y_i-x_i)\right|.
\end{equation}
By making use of \eqref{bounds}, \eqref{vb4} and \eqref{a1} we find that the second term of the r.h.s. of (\ref{new1}) is bounded above by $mnC_0(n+1)C_2\,\diam\, U_0\, Q^{v_1}$. In view of \eqref{U_0} and \eqref{C_2} the latter is no larger than $\tfrac12Q^{v_1}$. On the other hand, by \eqref{vb2} the first term in the r.h.s. of (\ref{new1}) is $\ge Q^{v_1}$. Thus, (\ref{new1}) implies that
$$
\left|\tfrac{\partial}{\partial x_1}(F+\theta)(\vv y)\right|\ge
\tfrac12Q^{v_1}\,.
$$
On the other hand, by  using \eqref{bounds}, \eqref{vb4} and \eqref{a1} we find that
$$
\left|\tfrac{\partial}{\partial x_i}(F+\theta)(\vv y)\right| \ \le\
C_4 \, Q^{v_1}
$$
for any $i=1,\dots,m$ and $\vv y\in U_0$, where $$C_4:=(n+1)C_0\max\{C_3,\,(n+1)C_2\}\, . $$
This together with the above lower bound inequality implies (\ref{extra}) with
$p:=(2mC_4)^{-1}$.

\bigskip

\noindent\textbf{Verifying that $\vv x\in\Delta(R_F,\rho(Q))$:} We will  makes use of the following easy consequence of the Mean Value Theorem.

\begin{lemma}\label{lem5+}
Let $f:I\to\R$ be a $C^1$ function on an interval $I$ such that $|f'(x)|\ge d>0$ for all $x\in I$. Let $x_1\in I$ and suppose that $B(x_1,|f(x_1)|d^{-1})\subset I$. Then, there is an $x_0\in B(x_1,|f(x_1)|d^{-1})$ such that $f(x_0)=0$.
\end{lemma}

\noindent Let $\vv x=(x_1,\dots,x_m)$. Consider the interval
$$
I:=\{x\in\R:(x,x_2,\dots,x_m)\in B\}
$$
and the function $f:I\to\R$ given by $f(x)=(F+\theta)(x,x_2,\dots,x_m)$.
In view of  (\ref{vb1}) and (\ref{vb2}) and the fact that  $\vv x\in\tfrac12B$, Lemma~\ref{lem5+} is applicable and implies that there exists some $x_0\in I$ such that $f(x_0)=0$ and $|x_1-x_0|\le (n+1)C_2 Q^{-n-v_1}$. Then $\vv x':=(x_0,x_2,\dots,x_m)\in B$ satisfies $F(\vv x')+\theta(\vv x')=0$ and
\begin{equation}\label{e:011}
|\vv x-\vv x'|\le (n+1)C_2 Q^{-n-v_1}\,.
\end{equation}

\noindent
On making use of (\ref{bounds}) and the Mean Value Theorem, we find that $|(F+\theta)(\vv y)|\ll Q^{-n}$  for any $\vv y$ satisfying  $|\vv y-\vv x'|\ll Q^{-n-v_1}$.  Then,  on using the above argument  for determining  $\vv x'$,  enables us to conclude  that for sufficiently large $Q$ the ball of radius $3\rho(\beta_F)$ centred at $\pi\vv x'$ is contained in $\tilde V$, where  $\pi$ is the projection map  given  by (\ref{pi}) and  $\tilde V$ is  as in (\ref{ressets}). The details are pretty straightforward and are  left to the reader. The upshot is  that $\vv x'\in R_F$ which together with  (\ref{e:011}) implies that
$\vv x\in \Delta(R_F,\rho(Q))$ as required, where $$
\rho(Q)=\kappa_1 Q^{-n-v_1}\qquad\text{with }\quad\kappa_1:=(n+1)C_2\,.
$$

\bigskip
\bigskip

\noindent\emph{Acknowledgements.} The authors are extremely grateful to Bob Vaughan for all his support over the years and for promoting Number Theory at York. Also happy number sixty five!

{\small

\def\cprime{$'$}
  \def\polhk#1{\setbox0=\hbox{#1}{\ooalign{\hidewidth
  \lower1.5ex\hbox{`}\hidewidth\crcr\unhbox0}}}

\noindent Dzmitry Badziahin: Department of Mathematics,
University of York,\\
\phantom{Dzmitry Badziahin: }Heslington, York, YO10 5DD,
England\\
\phantom{Dzmitry Badziahin: }\textit{E-mail}: \verb|db528@york.ac.uk|

\vspace{5mm}

\noindent Victor Beresnevich: Department of Mathematics, University of York,\\
\phantom{Victor Beresnevich: }Heslington, York, YO10 5DD,
England\\
\phantom{Victor Beresnevich: }\textit{E-mail}: \verb|vb8@york.ac.uk|

\vspace{5mm}

\noindent Sanju  Velani: Department of Mathematics, University of
York,\\
\phantom{Sanju  Velani: }Heslington, York, YO10 5DD,
England\\
\phantom{Sanju  Velani: }\textit{E-mail}: \verb|slv3@york.ac.uk|

\end{document}